\documentclass{lmcs}
\pdfoutput=1

\usepackage{lastpage}
\lmcsdoi{15}{1}{23}
\lmcsheading{}{\pageref{LastPage}}{}{}%
{Jul.~20,~2016}{Mar.~05,~2019}{}

\usepackage[utf8]{inputenc}
\usepackage{latexsym,amssymb}
\usepackage{graphicx}

\theoremstyle{plain}
\newtheorem{lemma}[thm]{Lemma}
\newtheorem{proposition}[thm]{Proposition}
\newtheorem{theorem}[thm]{Theorem}
\newtheorem{corollary}[thm]{Corollary}
\newtheorem{example}[thm]{Example}

\def\BP{\noindent {\it Proof.} }
\def\EP{\hfill$\Box$

}
\def\BC{\begin{center}}
\def\EC{\end{center}}
\def\BT{\begin{theorem}}
\def\ET{\end{theorem}}
\def\BPR{\begin{proposition}}
\def\EPR{\end{proposition}}
\def\BL{\begin{lemma}}
\def\EL{\end{lemma}}
\def\BCO{\begin{corollary}}
\def\ECO{\end{corollary}}
\def\BPI{\BC \begin{picture}}
\def\EPI{\end{picture}\EC}
\def\BTA{\begin{tabular}}
\def\ETA{\end{tabular}}
\def\BX{\begin{example}}
\def\EX{\end{example}}
\def\BIT{\begin{enumerate}[align=left]}
\def\EIT{\end{enumerate}}
\def\IT{\item}

\def\B{{\mathcal B}}
\def\C{{\mathcal C}}
\def\D{{\mathcal D}}

\def\F{{\mathcal F}}

\def\I{{\mathcal I}}

\def\O{{\mathcal O}}
\def\P{{\mathcal P}}

\def\S{{\mathcal S}}
\def\T{{\mathcal T}}
\def\U{{\mathcal U}}

\def\Y{{\mathcal Y}}
\def\Z{{\mathcal Z}}

\def\Ro{\hspace{-.2ex}\varrho}

\begin{document}
\author{Marcel Ern\'{e}}
\address{Leibniz University, Faculty for Mathematics and Physics, Hannover, Germany}
\email{erne@math.uni-hannover.de}
\urladdr{http://www2.iazd.uni-hannover.de/\textasciitilde erne/}
\title[Web spaces and worldwide web spaces]
{Web spaces and worldwide web spaces:\\ topological aspects of domain theory}
\begin{abstract}
Web spaces, wide web spaces and worldwide web spaces ({\it alias} C-spaces) provide useful generalizations of continuous domains. 
We present new characterizations of such spaces and their patch spaces, obtained by joining the original topology with a second topology having the dual specialization order; these patch spaces 
possess good convexity and separation properties and determine the original web spaces.  
The category of C-spaces is concretely isomorphic to the category of fan spaces; these are certain
quasi-ordered spaces having neighborhood bases of fans, where a fan is obtained by deleting a finite number of principal dual ideals from a principal dual ideal. 
Our approach has useful consequences for domain theory, because the T$_0$ web spaces are exactly the generalized Scott spaces associated with locally approximating ideal extensions, 
and the T$_0$ C-spaces are exactly the generalized Scott spaces associated with globally approximating and interpolating ideal extensions. 
The characterization of continuous lattices as meet-continuous lattices with T$_2$ Lawson topology and the Fundamental Theorem of Compact Semilattices are extended to non-complete posets.
Finally, cardinal invariants like density and weight of the involved objects are investigated. 
\end{abstract}

\maketitle

\section*{Introduction}
\label{intro}
Suitable models for the theory of computation and approximation are certain ordered or quasi-ordered sets ({\em qosets} for short), 
whose elements represent states of computation, knowledge or information, while the order relation abstractly describes refinement, improvement or temporal sequence.
Domain theory in its widest sense may be regarded as the order-theoretical framework for approximating relations. In the primary source for the present study, the monograph 
{\em Continuous Lattices and Domains} \cite{CLD}, an {\em approximating auxiliary relation} for a (partially) ordered set or {\em poset} $(X,\leq)$ is defined to be a relation $\prec$ on $X$ such that  
\BIT
\IT[] $w\leq x \prec  y \leq z \, \Rightarrow \, w \prec z \, \Rightarrow \, w\leq z$,
\EIT
and for each $y\in X$, the set ${\prec\!y} = \{ x \!: x \prec y\}$ is an ideal (i.e.\ a directed lower set) with join $y$. Such a relation is rarely reflexive but always transitive. Hence, if it has the {\em interpolation property}
\BIT
\IT[] $x \prec z \ \Rightarrow \exists\, y \ (x \prec y \prec z)$
\EIT
then it is idempotent. A study of such relations and their topological manifestation was presented in \cite{EABC}, 
where an idempotent relation $\Ro$ on a set $X$ was called a {\em $C$-quasi-order} if each set $\Ro y = \{ x : x\ \Ro\ y\}$ is an ideal with respect to the quasi-order $\leq_{\Ro}$ 
given by $x\leq_{\Ro} y \, \Leftrightarrow \, \Ro x \subseteq \Ro y$, and a {\em C-order} if, in addition, it has the {\em separation\,\,property} 
\BIT
\IT[] $\Ro x = \Ro y \, \Rightarrow \, x = y$, \ i.e., $\leq_{\Ro}$ is an (antisymmetric) order.
\EIT
The C-quasi-ordered sets are just the domain-theoretical {\em abstract bases} \cite{AJ},\,\cite{CLD},\,\cite{Lri}, \cite{Sm}, i.e.\ sets $X$ 
with a relation $\Ro$ on $X$ satisfying for all finite $F \subseteq X$ the equivalence
\BIT
\IT[] $F \subseteq {\Ro z}  \ \Leftrightarrow \ \exists \,y \in \Ro z \ (F\subseteq {\Ro y})$.  
\EIT
It is easy to check that each approximating auxiliary relation with the interpolation property is a C-order, and conversely, 
that each C-order $\Ro$ on a set $X$ is an interpolating and approximating auxiliary relation for the poset $(X,\leq_{\Ro})$. 
Thus, the order relation is important but in some sense redundant in the theory of approximation. 

In the mathematical and in the computer-theoretically oriented literature, the word {\em `domain'} represents quite diverse structures, 
and its meaning ranges from rather general notions like {\em dcpos} (up-complete posets), in which all directed subsets have suprema, 
to quite specific kinds of ordered sets like $\omega$-algebraic dcpos, sometimes with additional properties (see, e.g., \cite{AJ},\,\cite{Edom},\,\cite{CLD},\,\cite{SLG},\,\cite{V}). 
But even requiring the existence of suprema for all directed sets is often too restrictive for certain desired applications.
Therefore, as in \cite{CLD}, we adopt the convention to speak of a {\em continuous poset} if for each element $x$ there is a least ideal having a join above\,\,$x$ 
(the {\em way-below ideal} ${\ll\!y}$), but such a poset need not be a dcpo.
The {\em continuous domains} are then the continuous dcpos; observe that in \cite{Com} and \cite{Hoff2} they are called {\em continuous posets}, whereas in \cite{CLD} they are simply referred to as {\em domains}.

Although domains are usually defined in order-theoretical terms, there exist also topological descriptions of them. In \cite{EABC}, a concrete isomorphism has been established between the category of 
C-quasi-ordered sets (with suitable morphisms) and the category of {\em C-spaces}, i.e.\ spaces whose topology is completely distributive, generalizing the isomorphism between the category of continuous domains 
(resp.~their way-below relations) and that of sober C-spaces \cite{EScon},\,\cite{Hoff2},\,\cite{CLD},\,\cite{Law}.  
(Caution:\,in \cite{Kue} and \cite{Pri}, the word `{\em C-space}' has a different meaning.)
By a space we always mean a topological space, but extensions to arbitrary closure spaces are possible (see \cite{EABC},\,\cite{Epol},\,\cite{Eclo}). 
The lattice of closed sets of a space $(X,\S)$ with topology $\S$ is denoted by $\S ^c$, the closure of a subset $Y$ by $cl_{\S}Y$ or $Y^-\!$, and the interior by $int_{\S}Y\!$ or $Y^{\circ}$. 
A bridge between order and topology is built by the {\em specialization order\,}:
$$
x\leq y \ \Leftrightarrow \ x\leq_{\S} y \ \Leftrightarrow \ x\in \{ y\}^- \ \Leftrightarrow \  \forall\,U\! \in \S\ (x\in U \,\Rightarrow \, y\in U).
$$  
It is antisymmetric, hence an order, iff $(X,\S)$ is T$_0$, but we speak of specialization orders also in the non--T$_0$ setting.
The {\em saturation} of a subset $Y$ is the intersection of all its neighborhoods; this is the {\em up-closure} of $Y$, the upper set ${\uparrow\!Y}$ generated by $Y$ relative to the specialization order.
Accordingly, in the {\em specialization qoset}  $\Sigma^{-\!} (X,\S) = (X,\leq_{\S})$, the upper sets are the saturated sets, and the principal filters ${\uparrow\!x} = \{ y : x\leq y\}$ are the {\em cores}, 
that is, the saturations of singletons; on the other hand, the principal ideals ${\downarrow\!y} = \{ x : x\leq y\}$ are the point closures, and the lower sets $Y$ 
(satisfying $Y = {\downarrow\!Y} = \bigcup \,\{ {\downarrow\!y}:y\in Y\}$) are the unions of closed sets.

Several classes of spaces that interest us here may be characterized by certain infinite distributive laws for their lattices of open or closed sets \cite{Eweb}. 
Recall that a {\em frame} or {\em locale} (see, e.g., \cite{Jo} or \cite{V}) is a complete lattice $L$ satisfying the identity
\vspace{-.5ex}
$$
\textstyle{ {\rm (d)} \ x\wedge \bigvee Y = \bigvee \{ x\wedge y : y\in Y \} \vspace{-.5ex}}
$$
for all $x\in L$ and $Y \subseteq L$; the dual of (d) characterizes {\em coframes}. The identity
\vspace{-.5ex}
$$
\textstyle{{\rm (D)} \ \bigwedge \,\{ \bigvee Y : Y \!\in \Y\} = \bigvee\bigcap \Y} \vspace{-.5ex}
$$
for all collections $\Y$ of lower sets, defining {\em complete distributivity}, is much stronger.
However, frames may also be defined by the identity (D) for all {\em finite} collections $\Y$ of lower sets. 
An up-complete meet-semilattice satisfying (d) for all ideals (or directed sets) $Y$ is called {\em meet-continuous}. 
Similarly, the {\em continuous lattices} in the sense of Scott \cite{Com},\,\cite{Sco} are those complete lattices which enjoy the identity (D) for all collections of ideals.
Therefore, completely distributive lattices are also called {\em super\-continuous};
alternative descriptions of complete distributivity by equations involving choice functions are equivalent to {\sf AC}, the Axiom of Choice \cite[4.3--E\,10]{Herr}.

A complete lattice satisfying (D) for all collections of lower sets generated by finite sets is called {\em $\F$-distributive} \cite{EZcon} or a {\em wide coframe} \cite{Eweb},
and its dual a {\em quasi-topology} \cite{KP} or a {\em wide frame} \cite{Eweb}; the quasi-topologies are exactly the complete homomorphic images of topologies, while
the isomorphic copies of topologies are the {\em spatial frames}, in which every element is a meet of primes. Thus, spatial frames are wide frames, but not conversely (see K\v{r}i\v{z} and Pultr \cite{KP} for a counterexample). 

In \cite{Eweb}, we have introduced three classes of spaces that might be useful for the mathematical foundation of computation, communication and information theory: 
\begin{enumerate}[align=left]
\IT[--] in a {\em web space} each point has a neighborhood base of webs, i.e.\ unions of cores each of which contains the point,\\[-3ex]
\IT[--] in a {\em wide web space} each neighborhood $U$ of a point contains a neighborhood $V$ such that $V\!\dashv U$, i.e., each finite subset of $V$ has a lower bound in $U$,\\[-3ex]
\IT[--] in a {\em worldwide web space} each point has a neighborhood base of cores.
\end{enumerate}

\begin{figure}[h!] 
  \smallskip
\includegraphics[height=7.4cm]{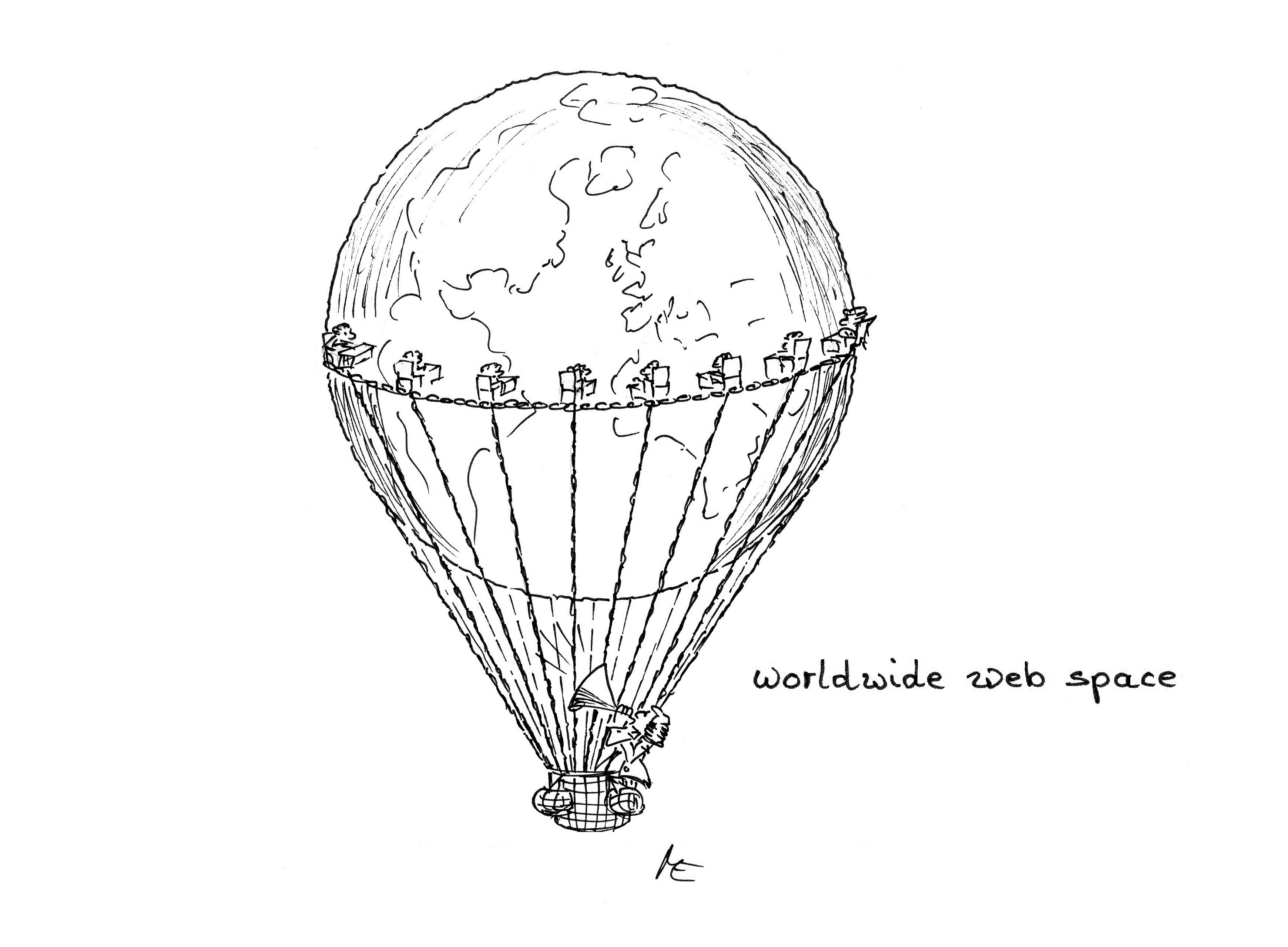}
\end{figure}
The following point-free characterizations are established in \cite{Eweb}: a space is a
\BIT
\IT[--] web space iff its topology is a coframe,
\IT[--] wide web space iff its topology is a wide coframe,
\IT[--] worldwide web space iff its topology is completely distributive.
\EIT
Many characterizations of web spaces are given in \cite{Eweb} and \cite{Epatch}. Some more basic facts and motivation to consider (wide or worldwide) web spaces will be given in Section\,\,\ref{basic}. 
Since spaces typically arising in domain theory often do not satisfy higher separation axioms (e.g., T$_1$ web spaces must already be discrete), passing to patch spaces with better separation properties appears fruitful. 
The patch spaces of a given space are obtained by joining its topology with a {\em cotopology} 
(or {\em complementary topology} \cite{Lss}), that is, a topology having the dual specialization order, and adding the original specialization order.
In Section \ref{convex}, we briefly review some applications of the construction of general patch spaces to web spaces, initiated in \cite{Epatch}.

Originally, worldwide web spaces have been termed {\em C-spaces} \cite{EScon},\,\cite{EABC}; the letter $C$ refers not only to {\bf c}ores and to {\bf c}omplete distributivity, 
but also to the facts that
\BIT
\IT[--] they are the locally {\bf c}ompact wide web spaces,
\IT[--] they are those spaces whose lattice of {\bf c}losed sets is {\bf c}ontinuous,
\IT[--] they are those spaces whose {\bf c}losure operator induces a {\bf c}omplete homomorphism from the {\bf c}omplete lattice of lower sets to that of {\bf c}losed sets.
\EIT
C-spaces share useful properties with the more restricted {\em B-spaces}, which have a least base, consisting of all open cores, and generalize algebraic domains, which correspond to the sober B-spaces. Still more limited 
are the {\em A-spaces} or {\em Alexandroff-discrete spaces} \cite{Alex}, in which all cores are open; the sober A-spaces represent the Noetherian domains. See the {\em ABC of Order and Topology} \cite{EABC} for details. 
In contrast to A- and B-spaces, the C-spaces are general enough to cover important examples of analysis: the real line ${\mathbb R}$ and its finite powers are continuous posets (but not dcpos\,!), hence C-spaces 
with the Scott topology, whose weak patch topology is the Euclidean topology. 
Section \ref{CDS} is devoted to a survey of the main properties of worldwide web spaces or C-spaces. 
Such spaces occur in diverse fields of mathematics\,--\,not only topological but also algebraic and order-theoretical ones\,--\,and in theoretical computer science. 
In particular, we present some new aspects of the categorical equivalence beetween C-spaces and C-quasi-ordered sets (abstract bases).

In Section \ref{sectorspaces}, we characterize the patch spaces of C-spaces as {\em sector spaces}. 
These are {\em ${\uparrow}$-stable semi-qospaces} (i.e.\ quasi-ordered spaces in which up-closures of open sets are open, and principal ideals and filters are closed) 
whose points have neighborhood bases of {\em sectors}, a special kind of webs possessing least elements. Restriction of the patch functors ${\rm P}_{\zeta}$ (see Section \ref{convex}) to the category of C-spaces yields 
concrete isomorphisms to the categories of so-called {\em $\zeta$-sector spaces}. 
In particular, the weak patch functor ${\rm P}_{\upsilon}$ induces a \mbox{categorical} isomorphism between C-spaces and {\em fan spaces},
i.e.\ ${\uparrow}$-stable semi-qospaces in which each point has a neighborhood base of {\em fans} ${\uparrow\!u}\setminus{\uparrow\!F}$ with finite sets $F$. 
In Section \ref{fanspaces}, such fan spaces are investigated and characterized by diverse order-topological properties, including uniformizability. 

The continuous domains, equipped with the Lawson topology, are just those fan spaces whose upper spaces are sober. In Section \ref{domain},
we find alternative descriptions of such ordered spaces, comprising a strong hyperconvexity property, high order-separation axioms and conditions on their interior operator. 
This enables us to generalize the characterization of continuous lattices as meet-continuous lattices whose Lawson topology is Hausdorff \cite[III--2.11]{CLD} to non-complete situations.

In Section \ref{semi}, we draw some new conclusions about semitopological or topological semilattices, and we extend the Fundamental Theorem of Compact Semilattices \cite[VI--3]{CLD} to the setting of posets.
Crucial is the fact that a semilattice $S$ with a topology having $S$ as specialization poset is semitopological iff it is a web space, and is a topological semilattice with small semilattices iff it is a wide web space.

In Section \ref{corebasis}, we note that the categories of T$_0$ C-spaces and of ordered fan spaces are not only isomorphic to the category of C-ordered sets, 
but also equivalent to the category of {\em based domains}, i.e.\ pairs consisting of a continuous domain and a basis of it (in the sense of \cite[III--4]{CLD}). 
Thus, the theory of C-spaces may be regarded as a part of domain theory, but the reverse containment makes sense as well. 

In Section \ref{power}, we revisit powerdomain constructions in the setting of C-spaces.

Finally, in Section \ref{density},\,we study weight and density of the spaces under consideration, using the order-theoretical description of C-spaces. 
For example, the weight of a C-space is equal to the density of its strong patch space, 
but also to the weight of the lattice of closed sets. This leads to the conclusion that the weight of a completely distributive lattice is always equal to the weight of the dual lattice.

\vspace{.3ex}

\noindent {\bf Acknowledgement.}  I am indebted to an anonymous referee who has read very carefully the first draft of this paper and has made many valuable suggestions that helped to improve the presentation considerably.

\section{Web spaces, wide web spaces and worldwide web spaces}
\label{basic}

We shall work in {\sf ZF} or {\sf NBG} (Zermelo--Fraenkel or Neumann--Bernays--G\"odel set theory) and point out where {\sf AC} or a weaker choice principle is needed.
For categorical notions see \cite{AHS}, for order-theoretical and topological definitions and facts refer to \cite{CLD}, and for additional material on domain theory to\,\cite{AJ}.

Recall that a quasi-ordered set or qoset is a pair $Q = (X\leq)$ with a reflexive and transitive relation $\leq$; its dual is denoted by $\geq$, and the dual qoset $(X,\geq)$ by\,\,$\widetilde{Q}$.
A {\em topological selection} $\zeta$ for qosets assigns to each qoset $Q$ a topology $\zeta Q$ whose specialization qoset is $Q$. 
Thus, $\zeta$ may be regarded as a functor right inverse to the specialization functor $\Sigma^-$.
The largest topological selection is $\alpha$, where $\alpha Q$ is the {\em Alexandroff topology} of all upper sets; the smallest one is $\upsilon$, 
where $\upsilon Q$ is the {\em weak upper topology} generated by the complements of principal ideals; and an intermediate one is $\sigma$, where $\sigma Q$ is the {\em Scott topology}, 
consisting of all upper sets $U$ that meet every directed subset $D$ of $Q$ having a least upper bound $y\in U$ (i.e.\ $D\subseteq {\downarrow\! z} \Leftrightarrow y\leq z$). 
The {\em weak patch topology} of $\sigma Q$, generated by $\sigma Q \cup \upsilon \widetilde{Q}$, is the {\em Lawson topology} $\lambda Q$.
We denote by $\Upsilon Q$ the {\em weak space} $(X,\upsilon Q)$, by $\Sigma Q$ the {\em Scott space} $(X,\sigma Q )$, and by $\Lambda Q$ the (quasi-ordered!) {\em Lawson space} $(Q,\lambda Q )$.

While web spaces and worldwide web spaces are quite natural and handy, the notion of wide web space is a bit subtle. A nonempty subset of a space is {\em strongly connected} 
if it is not contained in the union of two open sets unless it is already contained in one of them; in terms of the specialization order, that condition simply means that the set is filtered, i.e.\,down-directed.
A space is {\em locally strongly connected} if each of its points has a base of strongly connected neighborhoods, while a space {\em has a dual} \cite{Hoff1}, that is, its topology is dually isomorphic to another topology, 
if it has a base of (open) strongly connected sets. In the implication chain

{\it space with a dual $\,\Rightarrow $ locally strongly connected $\,\Rightarrow $ wide web space $\,\Rightarrow $ web space}

\noindent the first and the last implication cannot be inverted (see Example \ref{Ex11}), while it remains open whether the middle implication is an equivalence;
this is the case if the specialization qoset is a semilattice, or if the relation $\dashv $ is idempotent. 
Separating examples cannot be simple in view of the fact that `drawable' or `computable' posets must be countable,
if not finite, and that the following coincidences hold:

\BL
\label{countf}
{\rm (1)} Every countable frame is weakly atomic and spatial.\\
{\rm (2)} A countable frame that is also a coframe is algebraic and supercontinuous.
\EL    

\BP
(Cf.\ \cite{EDC}.) A poset is weakly atomic iff all its non-singleton intervals contain {\em jumps} $u \prec v$, meaning that $u < v$ holds but no $w$ satisfies $u < w < v$. Assume $L$ is a countable frame, 
but for some $x < y$ in $L$, the interval $[\,x,y\,] = \{ z : x\leq z\leq y\}$ contains no jump. Then one constructs recursively a countable dense subchain  of $[\,x,y\,]$, and by completeness of $L$,  
it would follow that $[\,x,y\,]$ contains a copy of the real unit interval, which is uncountable. Hence, $L$ must be weakly atomic, i.e., for $x < y$ there are $u,v$ with $x \leq u \prec v \leq y$.
The greatest element $p$ with $p \wedge v \leq u$ is then completely meet-irreducible,  hence prime, and $x\leq p$ but $y \not\leq p$.  This shows that $L$ is spatial.  
If $L$ is a coframe, too, then $p$ is even completely prime; the element $q = \bigwedge \{ r : r\not\leq p\}$ is then supercompact (i.e., $L\setminus {\uparrow\!q}$ is a principal ideal ${\downarrow\!p}$) with $q\leq y$ but
$q\not\leq x$. Hence, $L$ is superalgebraic, i.e.\ algebraic and supercontinuous.
\EP

Since the open cores are the supercompact members of topologies, we conclude:

\BCO
\label{contweb}
\mbox{A web space with countable topology is a B-space, hence a C-space.}
\ECO

We have seen that web spaces, wide web spaces and worldwide web spaces (i.e.\ C-spaces) may be characterized by certain infinite distributive laws 
that play a role in domain theory \,--\,but the converse is also true:

\BX
\label{Ex11}
{\rm Let $L = (X,\leq)$ be a distributive complete lattice and $\Upsilon L$ the weak space $(X,\upsilon L)$ with the complements of principal filters as subbasic open sets. Then}
\BIT
\IT[{\rm (1)}] $L$ is a frame iff $\Upsilon L$ is a web space,
\IT[{\rm (2)}] $L$ is a wide frame iff $\Upsilon L$ is a wide web space,
\IT[{\rm (3)}] $L$ is completely distributive iff $\Upsilon L$ is a worldwide web space.
\EIT
{\rm  We only prove (2). The cases (1) and (3) are similar but slightly simpler.

The following elementary remark will be helpful (cf.\ \cite{BaE}, \cite{EZcon}): a complete lattice is $\F$-distributive iff each of its elements $y$ is the join of all $x\dashv y$,
meaning $x\in {\downarrow\!F}$ for all finite $F$ with $y\leq \!\bigvee\! F$.
Dualizing that remark, we may restate the wide frame property by saying that for $x\not\leq y$ in $L$ there is a $z\in L$ with $x\not\leq z$ 
but $z\in {\uparrow\! G}$ whenever $G$ is a finite subset of $L$ with $\bigwedge G\leq y$. 
Now, if $x\in U = L\setminus {\downarrow\!F}$ for a finite $F\subseteq L$ then for each $y\in F$ we find a $z_y$ 
with the previous property; $E = \{ z_y : y\in F\}$ is then a finite subset of $L$ such that $x\in V = L\setminus {\downarrow\!E}$ and $\bigwedge G \in U$ for all finite $G\subseteq V$
(otherwise $\bigwedge G \leq y$ for some $y\in F$, hence $z_y\in {\uparrow\!G}\cap E$, in contrast to $G\cap {\downarrow\!E} = \emptyset$). Thus, these meets form a semilattice neighborhood of $x$ contained in $U$.
Hence, $\Upsilon L$ is locally strongly connected and so a wide web space.

Conversely, assume $\Upsilon L$ is a wide web space. For $x\not\leq y$ in $L$, the set $U = L\setminus {\downarrow\!y}$ is a $\upsilon$-open
neighborhood of $x$, so there is a basic $\upsilon$-open neighborhood $V = L \setminus {\downarrow\!F}$ of $x$ with $V\!\dashv U$. 
Then each finite set $G$ with $\bigwedge G\leq y$ intersects ${\downarrow\!F}$ (since $L\setminus {\downarrow\!F} = V \dashv U = L\setminus {\downarrow\!y}$). 
In order to prove the wide frame property of $L$, we have to find a $z\in L$ with $x\not\leq z$ but $z\in {\uparrow\!G}$ whenever $G$ is finite and $\bigwedge G\leq y$. We claim that at least one $z\in F$ has that property. 
Assume, on the contrary, that for each $z\in F$ there would exist a finite $G_z$ with $\bigwedge G_z\leq y$ but $z\not\in {\uparrow\!G_z}$. 
Then, by distributivity of $L$, the finite set $G = \{ \bigvee \chi [F] : \chi\in P = \prod_{z\in F} G_z\}$ would satisfy 
$\bigwedge G = \bigvee \{ \bigwedge G_z : z\in F\}\leq y$, hence $G \cap\, {\downarrow\!F} \not =\emptyset \not = F\cap\, {\uparrow\!G}$; 
but $z\in F\cap\, {\uparrow\!G}$ implies $\chi (z)\leq \bigvee \chi [F]\leq z$ for some $\chi\in P$, in contrast to $\chi (z)\in G_z$ and $z\not\in {\uparrow\!G_z}$.

Now let $L$ be a complete Boolean lattice, hence a frame.  Then, by (1), $\Upsilon L$ is a web space, while by (2) and (3), we have the following equivalences:

$\hspace{1.5ex}\Upsilon L$ is a worldwide web space $\,\Leftrightarrow $ $\Upsilon L$ is a wide web space  

$\Leftrightarrow $ $L$ is completely distributive $\,\Leftrightarrow L\,$ is a wide frame $\Leftrightarrow$ $L$ is atomic.

\noindent On the other hand, for any wide frame (in particular, for any topology) $L$ that is not completely distributive, $\Upsilon L$ is a wide web space but not a worldwide web space.\\[1.5ex]
Concerning products, we have the following facts, based on the Axiom of Choice: 
}
\EX

\BPR
\label{product}
A nonempty product of spaces is a web space (a wide web space, locally strongly connected) iff all factors are web spaces (wide web spaces, locally strongly connected) and all but a finite number of them are strongly connected.
Similarly, a nonempty product of spaces is a B-space (resp.\ C-space) iff all factors are B-spaces (resp.\ C-spaces) and all but a finite number of them are supercompact.
\EPR

\BP
Consider a family $(S_i : i\in I)$ of spaces with nonempty product $S = \prod_{i\in I}S_i$.
Suppose first that all $S_i$ are web spaces and only finitely many of them are not strongly connected. For any basic open neighborhood $U = \prod_{i\in I} U_i$ of some point $x$ in $S$, 
there is a finite set $F$ of indices such that $U_i = S_i$ is strongly connected for all $i\in I\setminus F$. For each $i\in F$, we find a web $W_i$ and an open $V_i$ such that $x_i\in V_i \subseteq W_i$.
Putting $V_i = W_i = S_i$ for $i\in I\setminus F$, we obtain a product-open set $V = \prod_{i\in I}V_i$ and a web $W = \prod_{i\in I} W_i$ with $x\in V\subseteq W \subseteq U$. Hence, $S$ is a web space.

Conversely, suppose $S$ is a web space. Then each factor $S_i$ must be a web space, too, because the $i$-th projection $p_i$ from $S$ to $S_i$ is a continuous open surjection, 
and the continuous open image of any web space is again a web space ($p_i$ preserves the specialization order). 
Now assume $S_i$ is not strongly connected for an infinite set $J$ of indices $i\in I$. Using {\sf AC}, pick $x,y$ in $S$
such that for no $j\in J$ the coordinates $x_j$ and $y_j$ have a common lower bound. Since $x$ has a web neighborhood $W$ in $S$, there is a product-open set $V = \prod_{i\in I}V_i$ such that $x\in V\subseteq W$ and 
$F = \{ i\in I : V_i \neq S_i\}$ is finite. Pick a $j\in J\setminus F$ and define $z\in S$ by $z_j = y_j$ and $z_i = x_i$ for $i\neq j$. By the hypothesis $V_j = S_j$, we have $z\in V\subseteq W$, so there is a $w\in W$ with $w\leq x$ and $w\leq z$, contradicting the fact that $x_j$ and $z_j = y_j$ have no common lower bound.

The case of wide web spaces is similar; we omit the details. The cases of locally strongly connected spaces, B-spaces and C-spaces have been settled in \cite{BrE}; 
note that a space is supercompact iff it has a least element in the specialization order.
\EP

A major motivation to consider web spaces and worldwide web spaces in the realm of domain theory is provided by the {\em Polarity between Approximation and Distribution}, discussed extensively in \cite{Epol}.
Crucial for topological aspects are {\em ideal extensions} of a poset $P = (X,\leq)$, i.e.\ collections of ideals containing at least all principal ideals. 
Any such ideal extension $\Z$ gives rise to a {\em generalized Scott topology}
$$\sigma_{\Z} P = \{ U \subseteq P : \forall \, I\in \Z \ (U \cap I = \emptyset \Rightarrow U \cap \Delta I = \emptyset)\}.$$
Here $\Delta I = \bigcap \,\{ {\downarrow\!z} : I \subseteq {\downarrow\!z}\}$ is the cut closure of $I$; notice that $x\in \Delta I$ means $x\leq \bigvee I$ if $I$ has a join. 
We say $\Z$ is {\em locally approximating} if for each $I\in \Z$ and each $x\in \Delta I$, there is a $J\in \Z$ with $x = \bigvee J$ and $J\subseteq I$. Similarly, $\Z$ is {\em globally approximating} 
if for each $x\in X$ there is a $J_x\in \Z$ with $x = \bigvee J_x$ and $J_x\subseteq I$ for all $I\in \Z$ with $x\in \Delta I$; and $\Z$ is {\em strongly approximating} if, in addition, it has the interpolation property
$\,x\in J_z  \Leftrightarrow  \exists \, y \in J_z \, (x\in J_y).$  
The difference between local and global approximation lies in the position of the quantifiers. By definition, $\Z$ is (strongly) globally approximating iff 
$P$ is a {\em (strongly) $\Z$-precontinuous poset} in the sense of \cite{EZcon} (cf.\ \cite{BaE}).
For the special ideal extension $\Z = \D_{\vee}$, consisting of all ideals possessing a join, the topology $\sigma_{\D_{\vee}}P$ is the usual Scott topology $\sigma P$; if $\D_{\vee}$ is locally approximating then $P$ is
meet-continuous in the sense of \cite[III-2]{CLD} and \cite{Kou}, and the converse holds in case $P$ is a meet-semilattice; similarly, $\D_{\vee}$ is (strongly) globally approximating iff $P$ is a continuous poset. 
The considerations in \cite{Epol} yield:

\BT
\label{Zcon}
The T$_0$ web spaces are exactly the generalized Scott spaces associated with locally approximating ideal extensions, and the T$_0$ C-spaces are exactly 
the generalized Scott spaces associated with strongly approximating ideal extensions. 
\ET

Thus, web spaces and C-spaces provide an optimal order-topological framework for an abstract theory of approximation.
For union complete global standard extensions, the interpolation property is a consequence of the global approximation property (see \cite{BaE} and \cite{EZcon} for details). 
This applies not only to the extension by all ideals, leading to the $s_2$-continuous posets and weak Scott topologies $\sigma_2 P$ discussed in \cite{EScon}, 
but also to the extension by all ideals possessing joins, leading to the original theory of continuous posets \cite{EScon}, \cite{Hoff2}, \cite{Law}. 
Even more interesting for practical purposes is perhaps the union complete extension $\C_{\omega}$ by all ideals that are generated by ascending sequences, or {\em $\omega$-chains} (cf.\ \cite[2.2.5]{AJ}).
Often, $\omega$-chains fit better with computational problems, because sequences and the associated {\em sequential Scott topology} $\sigma_{\C_{\omega}}P$ are easily handled, 
while directed sets may be rather complicated from the computational point of view. And, last but not least, it suffices to assume the existence of suprema for all $\omega$-chains and a bottom element 
in order to assure that all functions preserving joins of $\omega$-chains have a least fixpoint\,--\,an important tool in computational domain theory. 

\BX
\label{Ex12}
{\rm
The countable subsets of a given set form a (strongly) $\C_\omega$-continuous poset $P$. Indeed, every $y\in P$ is the union (join) of the set $J_y$ of all finite subsets of $y$, and this is a lower set generated by the 
$\omega$-chain $\{ y_n \!: n\in \omega\}$, where $y = \{ x_n \!: n\in \omega\}$ and $y_n = \{ x_k \!: k < n\}$. 
Hence, the sequential Scott topology $\sigma_{\C_{\omega}}P$ makes $P$ a worldwide web space, and the convergence of ascending sequences is topological.
}
\EX

These thoughts show that the approach via web spaces and C-spaces provides applications that would be excluded by demanding the existence of directed joins.

\section{Patch spaces of web spaces}
\label{convex}

A {\em (quasi-)ordered space} is merely a (quasi-)\-ordered set equipped with a topology. In this elementary definition,  
no separation properties and no relationship between order and topology are required. 
A quasi-ordered space is a {\em lower semi-qospace} if all principal ideals are closed, 
an {\em upper semi-qospace} if all principal filters are closed, and a {\em semi-qospace} if both conditions hold; these conditions mean that the quasi\-order is {\em lower semiclosed}, 
{\em upper semiclosed} or {\em semiclosed}, respectively, in the sense of \cite[VI-1]{CLD}. 
A semi-qospace with a partial order relation is referred to as a {\em semi-pospace} or a {\em T$_1$-ordered space\,}. 
A {\em qospace} is a space equipped with a closed quasi-order $\leq$, regarded as a subset of the square of the space; equivalently, for $x\not\leq y$ there are open $U$, $V$ with $x \in U$, $y \in V$
and ${\uparrow\!U} \mathop{\cap} {\downarrow\! V} = \emptyset$. A {\em pospace} is a qospace with a partial order \cite{CLD}.
Slightly stronger is the definition of {\em T$_2$-ordered spaces}, requiring that for $x \not\leq y$ there are an open upper set containing $x$ and a disjoint open lower set containing $y$. 
(Caution: some authors call such spaces {\em strongly $T_2$-ordered} and mean by a {\em T$_2$-ordered space} a pospace; cf.\ \cite{Kue}, \cite{McC}).
Clearly, T$_1$-ordered spaces are T$_1$, and pospaces are T$_2$ (Hausdorff). A quasi-ordered space is said to be {\em upper regular} if for each open upper set $O$ containing a point $x$,
there is an open upper set $U$ and a closed upper set $B$ such that $x\in U \subseteq B \subseteq O$, 
or equivalently, for each closed lower set $A$ and each $x$ not in $A$, 
there is an open upper set $U$ and a disjoint open lower set $V$ with $x\in U$ and $A\subseteq V$. 
An upper regular T$_1$-ordered space is said to be {\em upper T$_3$-ordered}. 
Note the following irreversible implications:

\bigskip
\begin{minipage}{.9\textwidth} 
{\em compact pospace $\Rightarrow$ upper\,T$_3$-ordered $\Rightarrow$\,T$_2$-ordered $\Rightarrow$ pospace $\Rightarrow$ semi-pospace\,}

\vspace{-0.5ex}
\hspace{8ex}$\Downarrow $\hspace{16ex}$\Downarrow $\hspace{31ex}$\Downarrow $\hspace{12ex}$\Downarrow $
\vspace{-0.5ex}

\indent {\em compact qospace $\,\Rightarrow$ upper regular semi-qospace $\hspace{4ex}\Rightarrow$ qospace\, $\Rightarrow$ semi-qospace\ }
\end{minipage}
\bigskip

\noindent For any quasi-ordered space $T = (X,\leq,\T)$ and the qoset $Q = (X,\leq)$,

\vspace{.5ex}

$\T^{\,\leq} = \T \cap \alpha Q\ $ is the topology of all open upper sets (in \cite{CLD} denoted by $\T^{\,\sharp}$),\\[1ex]
\indent $\T^{\,\geq} = \T \cap \alpha \widetilde{Q}\ $ is the topology of all open lower sets \hspace{.6ex}(in \cite{CLD} denoted by $\T^{\,\flat}$).\\[1ex]
We call $\hspace{.2ex}{\rm U} T = (X,\T^{\leq})$ the {\em upper space} and ${\rm L} T = (X, \T^{\geq})$ the {\em lower space} associated with $T$.
A basic observation is that for lower semi-qospaces, the specialization order of $\T^{\leq}$ is $\leq$, 
and for upper semi-qospaces, the specialization order of $\T^{\geq}$ is $\geq$.

Recall that a subset of a qoset is {\em (order) convex} iff it is the intersection of an upper set and a lower set. 
A quasi-ordered space is {\em locally convex} if the convex open subsets form a base,
and {\em strongly convex} if its topology $\T$ is generated by $\T^{\leq} \cup \T^{\geq}$.

The following definition will be crucial: a quasi-ordered space $(Q,\T)$ is said to be {\em upwards stable} ({\em ${\uparrow}$-stable})
if it satisfies the following equivalent conditions:  
\vspace{-1ex}
\BIT
\IT[{\rm (u1)}] $O\in \T$ implies ${\uparrow\!O}\in \T$.
\IT[{\rm (u2)}] $\T^{\leq} = \{ {\uparrow\! O} : O\in \T\}$.
\IT[{\rm (u3)}] The interior of each upper set is an upper set: $int_{\T} Y\! = int_{\T^{\leq}}Y$ if $Y \!= {\uparrow\!Y}$.
\EIT

For example, every semitopological semilattice is ${\uparrow}$-stable (see Section \ref{semi}). 

A {\em web} around a point $x$ in a qoset is a subset containing $x$ and with each point $y$ a lower bound of $\{ x,y\}$.
Any union of down-directed sets that contain a common point $x$ is a web around $x$, and conversely.
By a {\em web (quasi-)ordered space} we mean an ${\uparrow}$-stable (quasi-)ordered space in which every point has a neighborhood base of webs around it.
In the case of a space equipped with its specialization order, this becomes the definition of a web space.
In \cite{Epatch}, many characteristic properties of web quasi-ordered spaces are given, and web spaces are characterized by the property 
that their open sets are exactly the up-closures ${\uparrow\!U}$ of the open sets $U$ in any of their patch spaces.
Note that the {\em meet-continuous dcpos} in the sense of \cite{CLD} and \cite{Kou} are just those dcpos whose Scott space is a web space (see \cite{Eweb} and \cite{Epatch}).

We call a quasi-ordered space {\em locally lower bounded} if it is $\uparrow$-stable and each neighborhood $U$ of a point contains a neighborhood $V \!\dashv U$ of that point, 
where $V \! \dashv U$ means that every finite subset of $V$ has a lower bound in $U$. If the given quasi-order is the specialization order, this amounts to the definition of a wide web space. 
Similarly, we call a quasi-ordered space {\em locally filtered} if it is $\uparrow$-stable and each point has a base of neighborhoods that are filtered, i.e.\ down-directed. 

\BPR
\label{webpatch}
The strongly convex semi-qospaces are exactly the patch spaces (of their upper spaces). More specifically:

\noindent {\rm (1)} Patch spaces of web spaces $S$ are web quasi-or\-dered with upper space $S$.
Conversely, strongly convex web quasi-ordered semi-qospaces are patch spaces of their upper spaces, and these are web spaces.

\noindent {\rm (2)} Patch spaces of wide web spaces $S$ are locally lower bounded with upper space $S$. 
Conversely, strongly convex locally lower bounded semi-qospaces are patch spaces of their upper spaces, and these are wide web spaces. 

\noindent {\rm (3)}  Patch spaces of locally strongly connected spaces $S$ are locally filtered with upper space $S$. Conversely, 
strongly convex locally filtered semi-qospaces are patch spaces of their upper spaces, and these are locally strongly connected. 
\EPR

\vspace{-.5ex}

\BP
The first statement and (1) were established in \cite{Epatch}. 

(2) Let $S = (X,\S )$ be a wide web space and $T = (X,\leq, \T )$ a patch space of\,\,$S$. 
By (1), $T$ is ${\uparrow}$-stable, and $S$ is the upper space of $T$, i.e., $\S = \T^{\leq}$.
Given $x \in O \in \T$, find $U \!\in \S$ and $V \!\in\! \T^{\geq}$ with $x \in U \cap V \subseteq O$, 
and a $W \!\dashv U$ such that $x\in W\in \S$. Then, we have $W\cap V \dashv U\cap V$ (since $V = {\downarrow\!V}$), {\em a fortiori} $W\cap V \dashv O$. Thus, $T$ is locally lower bounded.
Conversely, let $T = (X, \leq, \T )$ be a strongly convex, locally lower bounded semi-qospace.
Then $T$ is web quasi-ordered and by (1) a patch space of the web space $(X,\T^{\leq})$.
For $x \in U \in \T^{\leq}$, there is a $V \!\dashv U$ with $x \in V \in \T$.
Then $x \in {\uparrow\! V} \!\dashv U$ and ${\uparrow\! V} \in \T^{\leq}$, by ${\uparrow}$-stability; thus, $(X,\T^{\leq})$ is a wide web space. 

(3) is shown by very similar arguments, so that we may omit the proof.
\EP

\BCO
\label{webco}
The strongly convex web T$_1$-ordered spaces are exactly the patch spaces of T$_0$ web spaces.
Similarly, the strongly convex, locally lower bounded T$_1$-ordered spaces are exactly the patch spaces of T$_0$ wide web spaces, 
and the strongly convex, locally filtered T$_1$-ordered spaces are exactly the patch spaces of locally strongly connected T$_0$ spaces.

\ECO

Specific patch constructions arise from so-called {\em coselections} $\zeta$, which choose for any topology $\S$ a subbase $\zeta \S$ of a cotopology $\tau_{\zeta}\S$, 
that is, a \mbox{topology} whose specialization order is dual to that of $\S$. The topology $\S^{\zeta}$ generated by $\S \mathop{\cup} \zeta \S$ is then a patch topology, 
and the associated (quasi-ordered!) {\em $\zeta$-patch space} is ${\rm P}_{\zeta}(X,\S) = (X, \leq_{\S} ,\S^{\zeta})$. 
A quasi-ordered space $(Q,\T)$ is said to be {\em $\zeta$-convex} if $\T$ is generated by $\T^{\leq} \mathop{\cup} \zeta (\T^{\leq})$. 

A variety of examples is obtained as follows: any topological selection $\zeta$ induces a coselection by putting $\zeta \S = \zeta \widetilde {Q}$ for any space $(X,\S )$ with specialization qoset $Q$. 
Of particular interest is the least topological selection $\upsilon$, choosing the weak upper topology $\upsilon Q$ and leading to the weak patch topology $\S^{\upsilon}$.
In that case, $\upsilon$-convexity amounts to the notion of {\em hyperconvexity}, which means that the sets $U\setminus {\uparrow\!F}$ with $U\in \T^{\leq}$ and finite $F$ form a base. 
The name reminds of the connection with the upper topology $\upsilon P$ and domain-theoretical hypercontinuity \cite[VII--3]{CLD}, but it
has not much to do with the classical definition of hyperconvex metric spaces.
Observe that hyperconvexity implies $\zeta$-convexity, which implies strong convexity, and the latter implies local convexity, but not conversely; 
separating examples are given in \cite{Epatch}. It is easy to see that hyperconvex upper regular semi-qospaces are regular.

A space $(X,\S)$ is said to be {\em $\zeta$-determined} if $\,\S^{\,\zeta\,\leq} = \S$.
A map between spaces is called {\em $\zeta$-proper} if it is continuous and $\zeta$-patch continuous; 
and a map between quasi-ordered spaces is {\em lower semicontinuous} if preimages of closed lower sets are closed. 
In \cite{Epatch}, many examples are given, and the following facts are established:

\vspace{-.5ex}

\BL
\label{pat}
The patch functor ${\rm P}_{\zeta}$ associated with a coselection $\zeta$ induces a concrete functorial isomorphism 
between the category of $\zeta$-determined spaces with continuous (resp.\ $\zeta$-proper) maps and that
of $\zeta$-convex semi-qospaces with isotone lower semicontinuous (resp.\ continuous) maps;
the inverse isomorphism is induced by the concrete upper space functor ${\rm U}$, sending a semi-qospace $(X,\leq, \T)$ to $(X,\T^{\leq})$. 
\EL

Now, invoking Proposition \ref{webpatch}, we derive from Lemma \ref{pat} the following conclusions:

\BT
\label{patzeta}
For any coselection $\zeta$, the patch functor ${\rm P}_{\zeta}$ induces concrete categorical isomorphisms between  \vspace{-1ex}
\BIT
\IT[{\rm (1)}] web spaces and $\zeta$-convex, web quasi-ordered semi-qospaces,
\IT[{\rm (2)}] wide web spaces and $\zeta$-convex, locally lower bounded semi-qospaces,
\IT[{\rm (3)}] locally strongly connected spaces and $\zeta$-convex, locally filtered semi-qospaces.
\EIT
\ET

\section{C-spaces and C-qosets (abstract bases)}
\label{CDS}

A subset $C$ of a space is super\-compact iff every open cover of $C$ has a member that contains $C$, or equivalently, the saturation of $C$ is a core. 
Hence, the locally supercompact spaces (in which every point has a base of supercompact neighborhoods) are nothing but the C-spaces (having core neighborhood bases at each point). 
Under the assumption of the Ultrafilter Theorem (a consequence of {\sf AC}), many properties of locally super\-compact spaces are shared by the more general {\em locally hypercompact spaces},
where a subset is hypercompact iff its saturation is generated by a finite set \cite {Eweb}. The equivalences (1)$\,\Leftrightarrow \,$(2)$\,\Leftrightarrow \,$(3) in the next proposition, proven in \cite{EABC} and \cite{Eweb}, motivated the term `worldwide web spaces' for C-spaces. 

\BPR
\label{supertop}
For a  space $S = (X,\S)$, the following conditions are equivalent:
\vspace{-.5ex}
\BIT
\IT[{\rm (1)}] $S$ is a worldwide web space (a C-space).
\IT[{\rm (2)}] $S$ is a locally hypercompact web space.
\IT[{\rm (3)}] $S$ is a locally compact wide web space.
\IT[{\rm (4)}] $S$ is a locally supercompact space.
\EIT
Moreover, they are equivalent to each of the following statements:
\BIT
\IT[{\rm (5)}] The lattice of open sets is supercontinuous (completely distributive).
\IT[{\rm (6)}] The lattice of closed sets is supercontinuous.
\IT[{\rm (7)}] The lattice of closed sets is continuous.
\IT[{\rm (8)}] The interior operator preserves arbitrary unions of upper sets.
\IT[{\rm (9)}] The closure operator preserves arbitrary intersections of lower sets.
\EIT
\EPR

On account of (8) resp.\,(9), the interior resp.\,closure operator of a C-space induces a complete homomorphism from the completely distributive lattice of upper resp.\,lower sets onto the lattice of open resp.\,closed sets.

Computationally convenient is the fact that C-spaces are in bijective correspondence to C-quasi-orders (see \cite{EABC} and the introduction). 
Notice that these are not really quasi-orders but idempotent relations $\Ro$ on a set $X$ such that the sets
\vspace{-.5ex}
$$\Ro y = \{ x \in X : x \mathrel{\Ro} y\} \ \ (y\in X)
$$ 
are ideals with respect to the {\em lower quasi-order} $\leq_{\Ro}$ defined by
\vspace{-.5ex}
$$
x\leq_{\Ro} y \ \Leftrightarrow \Ro x \subseteq \Ro y,
$$
and such a relation $\Ro$ is a {\em C-order} if, in addition, $\Ro x = \Ro y$ implies $x = y\,$. The pair $(X,\Ro )$ is referred to as a {\em C-qoset} or a {\em C-ordered set}, respectively. 
Recall that the C-qosets are just the domain-theoretical abstract bases (see \cite{AJ},\,\cite{CLD}, and Section\,\,\ref{corebasis}). 

For an arbitrary relation $\Ro$ on a set $X$ and any subset $Y$ of $X$, put
$$\Ro\, Y = \{ x\in X : \exists\,y\in Y\, (x\mathrel{\Ro}y)\} , \  Y\Ro = \{ x\in X : \exists\,y\in Y\, (y\mathrel{\Ro}x)\},\vspace{-.5ex}
$$
$$\ _{\Ro}\O = \{ \,\Ro\, Y : Y\subseteq X\} ,  \ \O_{\Ro} = \{ Y\Ro : Y\subseteq X\}. \  
$$
The systems $_{\Ro}\O $ and  $\O_{\Ro}$ are closed under arbitrary unions,
and if the sets $\Ro y$ are ideals (relative to $\leq_{\Ro}$) then $\O_{\Ro}$ is even a topology.
On the other hand, if $\Ro$ is idempotent then $_{\Ro}\O$ consists of all rounded subsets, 
where a subset $Y$ is said to be {\em round(ed)}  if $Y = \Ro Y$ (see \cite{GK},\,\cite{KL},\,\cite{Lri}).
Typical examples of C-orders are the {\em way-below relations} $\ll$ of continuous posets $P$, 
in which for any element $y$ the set
$$\textstyle{\ll\! y = \{ x\in P : x \ll y\} = \bigcap \,\{ D\in \D_{\vee} P : y\leq \bigvee\! D\}}
$$
is an ideal (the {\em way-below ideal\,}) with join $y$; recall that $\D_{\vee} P$ denotes the set of all ideals possessing a join. 
The next straightforward result leans on the interpolation property of continuous posets, saying that their way-below relation is idempotent.

\BL
\label{conti}
For every continuous poset, the way-below relation is a C-order $\Ro$ such that for all directed subsets $D$ (relative to $\leq_{\Ro}$), 
$y = \bigvee\! D$ is equivalent to $\Ro y = \Ro D$; and any C-order with that property is the way-below relation of a continuous poset. 
\EL

Every space $S = (X,\S)$ carries a transitive {\em interior relation} $\Ro_{\S}$, given by
\vspace{-.5ex}
$$
x \,\Ro_{\S}\,y \, \Leftrightarrow \, y \in ({\uparrow\!x})^{\circ} = int_{\S} ({\uparrow\!x}), \vspace{-.5ex}
$$
where ${\uparrow\!x} = \bigcap\, \{ U \in \S : x \in U\}$ is the core of $x$. Note that for any subset $A$ of $X$,\\
$\Ro_{\S}A$ is a lower set in $(X, \leq_{\S})$: $\,x \leq_{\S} y\ \Ro_{\S} \, z$ entails  $z\in ({\uparrow\!y})^{\circ} \!\subseteq\! ({\uparrow\!x})^{\circ}$, hence $x \, \Ro_{\S} \,z$.\!

A map $f$ between C-qosets $(X,\Ro)$ and $(X',\Ro')$ {\em interpolates} if $\Ro ' f(y) \subseteq \Ro'f[\hspace{.2ex}\Ro y]$ for all $y \in X$; 
and $f$ is {\em isotone} if it preserves the lower quasi-orders, i.e., $x\leq_{\Ro} y$ implies $f(x) \leq_{\Ro '}\! f(y)$.
Note that $f$ is interpolating and isotone iff $\Ro' f(y) = \Ro'f[\hspace{.2ex}\Ro y]$ (cf.\ \cite{EABC}).

\BT
\label{Cspace}
{\rm (1)} $(X,\S)$ is a C-space iff there is a C-quasi-order $\Ro$ with $\S = \O_{\Ro}$.
In that case, $\Ro$ is the interior relation $\Ro_{\S}$, and $\leq_{\Ro}\!$ is the specialization order $\leq_\S$.

{\rm (2)} By passing from $(X,\S)$ to $(X,\Ro_{\S})$, and in the opposite direction from $(X,\Ro)$ to $(X,\O_{\Ro})$,
the category of C-spaces and continuous maps is concretely isomorphic to the category of C-qosets (alias abstract bases) and interpolating isotone maps.

{\rm (3)} For each closed set $A$ in a C-space $(X,\S)$ with interior relation $\Ro$, the set $\Ro A$ is the least lower set with closure $A$.
The closure operator induces an isomorphism between the super\-continuous lattice $_{\Ro}\O$ of all rounded sets and that of all closed sets, 
while $\O_{\Ro}$ is the dually isomorphic supercontinuous lattice of all open sets. 

{\rm (4)} The irreducible closed subsets of a C-space are exactly the closures of directed sets (in the specialization qoset). 
For any irreducible closed set $C$, the set $\Ro\, C$ is the least ideal with closure $C$. The 
closure operator induces an iso\-morphism between the continuous domain of rounded ideals and that of irreducible closed sets.

{\rm (5)} The topology of a C-space $(X,\S)$ is always finer than the Scott topology of its specialization qoset.

{\rm (6)} The cocompact topology $\tau_{\pi}\S$ of a C-space $(X,\S)$, generated by the coselection $\pi \S$ of all complements of compact saturated sets, 
is the weak lower topology of the specialization qoset. The patch topology $\S^{\pi} = \S \vee \tau_{\pi}\S$ is the weak patch topology\,\,$\S^{\upsilon}$.  
\ET

\BP
(1) and (2) have been established in \cite{EABC}. 
There are several reasonable alternative choices for the morphisms. For example, the isotone quasi-open maps 
(for which the saturations of images of open sets are open \cite{HM}) between C-spaces 
are the relation preserving isotone maps between the associated C-qosets (see \cite{EABC}). 

(3) Let $\Ro$ be the interior relation of $(X,\S)$.
We prove for $A \subseteq X$ the identity $(\Ro A)^-\! = A^-$, using idempotency of $\Ro$ in the equivalence * below:

$y\in (\Ro A)^- \ \Leftrightarrow \ \forall\, U\in \O_{\Ro}\ (y\in U \Rightarrow U\cap \Ro A \neq \emptyset)$

$\Leftrightarrow  \ \forall \, x\in X\ (x \mathrel{\Ro}y \,\Rightarrow\, x\Ro \cap \Ro A \neq \emptyset) \ \stackrel{*}{\Leftrightarrow} \ 
\forall \, x\in X\ (x \mathrel{\Ro}y \Rightarrow x\Ro \cap A \neq \emptyset)$

$\Leftrightarrow \ \forall\, U\in \O_{\Ro}\ (y\in U \Rightarrow U\cap A \neq \emptyset) \ \Leftrightarrow \ y\in A^-.$

\noindent In particular, $(\Ro A)^- = A$ in case $A$ is closed. On the other hand, 
any rounded set $Y$ satisfies the identity $Y = \Ro Y = \Ro (Y^-)$: since each $x\Ro$ is an open set, we have 

$x\in \Ro Y \ \Leftrightarrow \ x\Ro \cap Y \neq \emptyset \ \Leftrightarrow \ x\Ro \cap Y^- \neq \emptyset \ \Leftrightarrow \ x\in \Ro (Y^-).$ 

\noindent And if $Y$ is any lower set with $Y^-\! = A$ then $\Ro A = \Ro (Y^-) = \Ro Y \subseteq {\downarrow\!Y} = Y.$
That the map $Y\mapsto Y^-$ yields an isomorphism between $_{\Ro}\O$ and $\S^c$ is now straightforward.

(4) A subset $C$ of a space is irreducible iff it is nonempty and $C \subseteq A \cup B$ implies $C \subseteq A$ or $C \subseteq B$ for any closed sets $A,B$.
Directed subsets and their closures are irreducible. 
The rounded ideals of a C-qoset (or a C-space) form a dcpo $\I_{\,\Ro}$,
being closed under directed unions.  It is easy to see that $I \ll J$ holds in $\I_{\,\Ro}\!$ iff there is an $ x\in J$ with $I\subseteq {\Ro x}$,
and that $J$ is the join of the ideal $\ll \!\!J$ in $\I_{\,\Ro}$. Thus, $\I_{\,\Ro}$ is a continuous domain (cf.\ \cite{GK},\,\cite{KL},\,\cite{Lss}).
For the isomorphism claim, observe that the coprime members of $_{\Ro}\O$ are the rounded ideals, the coprime closed sets are the irreducible ones, 
and a lattice isomorphism preserves coprimeness.

(5) By (4), $I = \Ro y = \Ro \hspace{.2ex}{\downarrow\hspace{-.4ex}y}$ is an ideal with ${\downarrow\! y} = I^-$, so $y$ is a least upper bound of $I$; hence,
$y \in U \in \sigma (X,\leq)$ implies $U\cap \Ro y \not = \emptyset$, i.e.\ $y\in U\Ro$. Thus, $U = U\Ro \in \S$.

(6) Let $C$ be a compact saturated set and $y \!\in\! X\setminus C$. 
Then $ U = X\setminus {\downarrow\!y}$ is an open neighborhood of $C$ not containing $y$. As $(X,\S)$ is a C-space, we
have \mbox{$U \!= \bigcup\, \{ int_{\S}({\uparrow\! x}) : x\in U\}$}, and by compactness of $C$, we find a finite $F \subseteq U$ with
$C \!\subseteq \bigcup\, \{ int_{\S} ({\uparrow\! x}) : x \!\in\! F\} \!\subseteq\! {\uparrow\! F} \!\subseteq\! U$. 
For $Q = (X,\leq)$, the set ${\uparrow\! F}$ is $\upsilon\widetilde{Q}$-closed and excludes $y$. 
Thus, $C$ is $\upsilon\widetilde{Q}$-closed, and $\tau_{\pi}\S$ is contained in $\upsilon\widetilde{Q}$; the reverse inclusion is clear, since cores are compact, saturated, and form a closed subbase for $\upsilon\widetilde{Q}$. 
\EP

The arguments for (4) and (5) show that for any irreducible closed set $C$ in a C-space, the set $\Ro\, C$ is not only the least ideal but also the least irreducible lower set with closure $C$.

Recall that a {\em sober} space is a T$_0$ space whose only irreducible closed sets are the point closures; 
and a {\em d-space} \cite{Wy} is a T$_0$ space in which the closure of any directed subset is the closure of a singleton.
The d-spaces are the {\em monotone convergence spaces} in \cite{CLD}. 
From Proposition \ref{supertop} and Theorem \ref{Cspace}, we deduce (cf.\,\cite{EABC},\,\cite{Hoff2},\,\cite{Law},\,\cite{Lss}):

\BCO
\label{contdom}
The following conditions on a space $S$ are equivalent:
\vspace{-.5ex}
\BIT
\IT[{\rm (1)}] $S$ is a sober C-space.
\IT[{\rm (2)}] $S$ is a locally supercompact d-space.
\IT[{\rm (3)}] $S$ is a locally compact d-space and a wide web space.
\IT[{\rm (4)}] $S$ is the Scott space of a (unique) continuous domain. 
\EIT  
\vspace{-.5ex}
The category of sober C-spaces and the concretely isomorphic category of continuous domains are dual to the category of supercontinuous spatial frames.
\ECO

For the case of continuous posets that are not necessarily domains, see \cite{EScon}, \cite{Emin}.
Afficionados of domain theory might remark that continuous frames are automatically spatial (see \cite{CLD}, \cite{HL}). 
But that `automatism' requires choice; and the spatiality of supercontinuous frames also leans on a weak choice principle \cite{ECD}.

Since a T$_0$ space and its sobrification have isomorphic open set frames, it follows from Corollary \ref{contdom} that
a T$_0$ space is a C-space iff its sobrification is the Scott space of a continuous domain. This completion process is reflected, via Theorem \ref{Cspace},
by the fact that the rounded ideal completion of a C-ordered set is a continuous domain, and the C-order $\Ro$ 
extends to the completion, which means that $x\mathrel{\Ro}y$ is equivalent to $\Ro x \ll \Ro y$ in the completion (cf. \cite{Emin}, \cite{CLD}, \cite{GK}, \cite{Lss}, \cite{Lri}).

\section{C-stable spaces and sector spaces}
\label{sectorspaces}

As every C-space is a web space, it is equal to the upper spaces of its patch spaces and is therefore $\zeta$-determined for any coselection $\zeta$ (Proposition \ref{webpatch} and Theorem\,\,\ref{patzeta}).
We now are going to determine explicitly these patch spaces, which turn out to have very good separation properties
(whereas the only T$_1$ C-spaces are the discrete ones). For alternative characterizations of such patch spaces, we need further properties of their interior operators. 
Call a quasi-ordered space $(Q, \T )$ {\em C-stable} if
$$
\textstyle{{\uparrow\!O} = \bigcup \,\{ \,int_{\T}({\uparrow\!u}) : u\in O\}  \ \mbox{ for each }O \in \T, }
$$
and {\em d-stable} if for any filtered (i.e.\ down-directed) subset $D$ of $Q$,
$$
\textstyle{int_{\T}D \subseteq \bigcup\,\{ int_{\T} \, ({\uparrow\!u}): u\in cl_{\T^{\geq}} D\}.}
$$
While C-stability is a rather strong property, d-stability is a much weaker one, trivially fulfilled if all dual ideals are principal. The terminology is justified by

\BL
\label{cstable}
Let $T = (X,\leq, \T)$ be a semi-qospace.

{\rm (1)} $T$ is d-stable whenever its lower space is a d-space.

{\rm (2)} $T$ is C-stable iff it is $\uparrow$-stable and its upper space is a C-space.

{\rm (3)} $T$ is C-stable iff it is upper regular, locally filtered and d-stable.   
\EL

\BP
(1) If the lower space $\,{\rm L}T$ is a d-space then for any $D$ that is filtered in $(X,\leq)$, i.e.\ directed in $(X,\geq)$, there is a $u$ 
with ${\uparrow\!u} = cl_{\T^{\geq}} D$, hence $\,int_{\T} D \subseteq int_{\T} ({\uparrow\!u})$. 

(2) Clearly, a C-stable semi-qospace $(X,\leq ,\T)$ is ${\uparrow}$-stable. 
For $x\in O\in \T^{\leq}$ there exists a $u\in O$ with $x\in int_{\T} ({\uparrow\!u})$.
Then $U = int_{\T}({\uparrow\!u})\in \T^{\leq}$ (by ${\uparrow}$-stability) and $x\in U\subseteq {\uparrow\!u} \subseteq {\uparrow\!O} = O$.
Since $\leq$ is the specialization order of $\T^{\leq}$, this ensures that $(X,\T^{\leq})$ is a C-space.

Conversely, suppose $(X,\T^{\leq}) = (X, \{ {\uparrow\!O} : O \in {\mathcal T}\})$ is a C-space. 
Then, for $O\in {\mathcal T}$ and $y \in {\uparrow\!O}$, there is an $x\in {\uparrow\!O}$ and a $U\in \T^{\leq} \subseteq \T$ with $y\in U\subseteq {\uparrow\!x}$.
Now, pick a $u\in O$ with $x\in {\uparrow\!u}$; then $y\in U \subseteq {\uparrow\!x} \subseteq {\uparrow\!u}$. Thus, 
${\uparrow\!O} \subseteq \bigcup \,\{ \,int_{\T}({\uparrow\!u}) : u\in O\}$. 

(3) Let $T = (X,\leq, \T )$ be a C-stable semi-qospace. For $x\in O\in \T^{\leq}$, there is a $u\in O$ with $x\in U = int_{\T} ({\uparrow\!u}) \subseteq {\uparrow\!u} \subseteq O$.
By ${\uparrow}$-stability, we have $U \in \T^{\leq}$, and since $T$ is a semi-qospace, ${\uparrow\!u}$ is a closed upper set. Hence, $T$ is upper regular. 

In order to check local filteredness, pick for $x\in O\in \T$ an element $u\in O$ with $x\in int_{\T} ({\uparrow\!u}) \mathop{\cap} O \subseteq {\uparrow\!u} \mathop{\cap} O\,$; 
this is a filtered set, possessing the least element\,\,$u$.

Moreover, $T$ is d-stable, since for any subset $D$ and any $x\in O = int_{\T} D$, 
there is a $u \in O \subseteq D \subseteq cl_{\T^{\geq}} D$ with $x\in int_{\T}({\uparrow\!u})$.

Conversely, suppose that $T$ is an upper regular, locally filtered (in particular ${\uparrow}$-stable) and d-stable semi-qospace. 
Then, for $O\in \T$ and $x\in {\uparrow\!O}$, we have:

$x\in {\uparrow\!O} \in \T^{\leq}$ \hfill (by ${\uparrow}$-stability),

there are $U\in \T^{\leq}$ and $B\in \T^{\geq c}$ with $x\in U \subseteq B \subseteq {\uparrow\!O}$ \hfill (by upper regularity),

a filtered $D$ with $x\in int_{\T} D \subseteq D \subseteq U$ \hfill (by local filteredness),

and an element $y\in cl_{\T^{\geq}} D$ with $x \in int_{\T} ({\uparrow\!y})$ \hfill (by d-stability).

\noindent It follows that $y \in cl_{\T^{\geq}} D \subseteq cl_{\T^{\geq}} U \subseteq B \subseteq {\uparrow\!O}$.
Now, choose a $u\in O$ with $u\leq y$, hence ${\uparrow\!y} \subseteq {\uparrow\!u}$, to obtain $x\in int_{\T} ({\uparrow\!u})$. 
Thus, we see that ${\uparrow\!O}$ is contained in $\bigcup\,\{ int_{\T} ({\uparrow\!u}) : u\in O\}$, showing that $T$ is C-stable.
\EP

\vspace{1ex}

By Lemma \ref{cstable}, every C-stable semi-qospace is a qospace (being upper regular), and C-stability of a semi-qospace splits into the following four properties: \vspace{.5ex}

(c1) {\em upper regular}, (c2) {\em bases of filtered neighborhoods}, (c3) {\em ${\uparrow}$-stable}, (c4) {\em d-stable}.\vspace{.5ex}

\noindent These properties are independent: none of them follows from the other three.

\BX
\label{Ex41}
{\rm
Let $L$ be a non-continuous wide frame, for example, a T$_2$ topology that is not locally compact, like the Euclidean topology on the rational line $\mathbb Q$. 
The arguments in Example \ref{Ex11} show that $L$ with the weak upper topology $\upsilon L$ has small semilattices, i.e., each point has a neighborhood base of subsemilattices. 
From this observation, one easily deduces that the interval topology $\iota L = \upsilon L \vee \upsilon \widetilde{L}$ makes $L$ a locally filtered topological meet-semilattice ${\rm I} L$, hence ${\uparrow}$-stable
(cf.\ Section \ref{semi}); 
and ${\rm I}L$ is d-stable, since by the Ultrafilter Theorem, the lower space of ${\rm I}L$ is $\Upsilon \widetilde{L}$ (cf.\,\cite[III--3.18]{CLD}), which is a d-space for any complete lattice $L$. 
But the semi-pospace ${\rm I}L$ can be upper regular only if it is a pospace, hence T$_2$, which happens only if $L$ is a continuous lattice \cite[III--3]{CLD}. Thus, ${\rm I}L$ satisfies (c2), (c3) and (c4), but not (c1). 
}
\EX

\BX
\label{Ex42}
{\rm 
For $0 < s < 1$, consider the non-compact subspace $S = [\,0, s\,[ \, \mathop{\cup}\, \{ 1\}$ of the Euclidean space $({\mathbb R},\T)$, ordered by
$x \sqsubseteq y \Leftrightarrow x \!=\! y$ or $x \!=\! 0$ or $y \!=\! 1$. As $\{ 1\}$\linebreak 
is clopen, it is easy to see that $S$ satisfies (c1), (c3) and (c4) (all dual ideals are principal), 
but not\,(c2): no point except\,$1$\,has a filtered neighborhood not \mbox{containing\,0.} 
}
\EX

\BX
\label{Ex43}
{\rm
In Example 4.1 of \cite{Epatch} it is shown that $L = \{ a,\top\} \cup \{ b_n : n\in \omega\}$, ordered by 
$x\leq y \ \Leftrightarrow   x = y \mbox{ or } x = b_0 \mbox{ or } y = \!\top \mbox{ or } x = b_i, y = b_j, i < j,$
is a complete but not meet-continuous lattice and a compact pospace when equipped with the Lawson topology. 
It satisfies (c1), (c2), and (c4) (all dual ideals are principal).
However, (c3) is violated, since $\{ a \}$ is open, while ${\uparrow\!a} = \{ a, \top\}$ is not. 
}
\EX

\begin{picture}(300,63)(150,0)

\put(-70,0){
\begin{picture}(200,70)

\put(286,31){\circle{4}}

\put(288,32){\line(1,1){25.5}}
\put(288,30){\line(1,-1){25.5}}
\put(288.5,30.5){\line(1,0){25}}
\put(288.5,31){\line(1,0){25}}
\put(288.5,31.5){\line(1,0){25}}
\put(315,31){\circle{4}}

\put(315,37){\line(0,1){20}}
\put(315,5){\line(0,1){20}}

\put(315,59){\circle*{4}}
\put(319,59){$1$}
\put(315,3){\circle*{4}}
\put(319,-1){$0$}
\put(319,29){$s$}

\put(345,29){$S$}
\end{picture}
}

\put(70,0){
\begin{picture}(200,70)

\put(286,31){\circle*{4}}
\put(275,29){$a$}
\put(287.5,32){\line(1,1){26}}
\put(287.5,30){\line(1,-1){26}}

\put(315,36){\circle*{4}}
\put(315,20){\circle*{4}}
\put(315,34){\line(0,-1){12}}
\put(315,18){\line(0,-1){13}}

\put(315,39){\line(0,1){5}}
\put(315,47){\circle{1}}
\put(315,51){\circle{1}}
\put(315,55){\circle{1}}

\put(315,59){\circle*{4}}
\put(319,57){$\top$}
\put(315,3){\circle*{4}}
\put(320,31){$b_2$}
\put(320,14){$b_1$}
\put(320,-1){$b_0$}

\put(350,29){$L$}
\end{picture}
}

\end{picture}

\noindent 
\BX
\label{Ex44}
{\rm
In the lattice ${\mathbb R}^{\,\omega}$ of real sequences with the component\-wise order, the way-below relation is empty, 
whence the Scott space $\Sigma ({\mathbb R}^{\,\omega})$ is not a C-space.
The Lawson space $\Lambda  ({\mathbb R}^{\,\omega})$ heavily differs from the product space $(\Lambda {\mathbb R})^{\,\omega}$, 
although both spaces have the same upper space $\Sigma ({\mathbb R}^{\,\omega}) = (\Sigma \mathbb R)^{\,\omega}$ \cite{EScon}, 
which is locally strongly connected, hence a wide web space. 
Both spaces, $\Lambda ({\mathbb R}^{\,\omega})$ and $(\Lambda {\mathbb R})^{\,\omega}$, satisfy (c2) and (c3), possessing small subsemilattices and a continuous meet operation. 
While the Lawson space $\Lambda ({\mathbb R}^{\,\omega})$ is irreducible (has no disjoint nonempty open subsets), 
the product space is even upper T$_3$-ordered; so, by Lemma \ref{cstable}, it cannot fulfil (c4).
For {\em uncountable} $I$, the Scott space $\Sigma ({\mathbb R}^{\,I})$ differs from the product space $(\Sigma {\mathbb R})^{\,I}$ \cite{EScon}.

}
\EX

Given a quasi-ordered space $(X,\leq,\T)$, we call a nonempty set of the form ${\uparrow\!u}\mathop{\cap} V$ a {\em sector} if $V$ belongs to $\T^{\,\geq}$,
and a {\em $\zeta$-sector} if $V$ may be chosen in $\zeta (\T^{\leq})$, for a coselection $\zeta$.
Hence, $u$ is the least element of the sector, and every sector is obviously a filtered web around any point it contains.
By a ($\zeta$-){\em sector space} we mean an ${\uparrow}$-stable semi-qospace in which
each point $x$ has a base of ($\zeta$-)sector neighborhoods ${\uparrow\!u}\cap V$ (but the point $x$ need not be the minimum $u$ of such a sector).

\BPR
\label{sectors}

The sector spaces are exactly the
\vspace{-.5ex}
\BIT
\IT[{\rm (1)}] patch spaces of C-spaces,
\IT[{\rm (2)}] strongly convex, C-stable semi-qospaces,
\IT[{\rm (3)}] strongly convex, upper regular, locally filtered and d-stable qospaces. 
\EIT
In particular, all ordered sector spaces are upper T$_3$-ordered, a fortiori T$_2$-ordered.\\
Specifically, for any coselection $\zeta$, the $\zeta$-sector spaces are exactly the
\vspace{-.5ex}
\BIT
\IT[{\rm (1$\zeta$)}] $\zeta$-patch spaces of C-spaces,
\IT[{\rm (2$\zeta$)}] $\zeta$-convex, C-stable semi-qospaces,
\IT[{\rm (3$\zeta$)}] $\zeta$-convex, upper regular, locally filtered and d-stable qospaces.  
\EIT
\EPR

\BP
Sector spaces $T = (X,\leq,\T)$ satisfy (2): For $x\in O\in \T$, there are $u \!\in\! X$ and $V \!\in\! \T^{\,\geq}$ such that 
$x\in int_{\T} ({\uparrow\!u} \mathop{\cap} V) = int_{\T} ({\uparrow\!u})\cap V \subseteq {\uparrow\!u} \mathop{\cap} V \subseteq O$.
Then $u\in {\uparrow\!u}\cap V \subseteq O$ and $x \in int_{\T}({\uparrow\!u})$,
whence $O \subseteq \bigcup \,\{ int_{\T}({\uparrow\!u}) : u \!\in\! O\}$.
Using ${\uparrow}$-stability and applying the previous argument to ${\uparrow\!O}$ instead of $O$, one concludes that $T$ is C-stable.
Again by ${\uparrow}$-stability, $U = {\uparrow int_{\T}}({\uparrow\!u}\,\cap V)$ belongs to $\T^{\,\leq}$, and the above reasoning yields 
$x\in U\cap V \subseteq {\uparrow\!u}\cap V\subseteq O$, proving strong convexity.

(2)$\,\Rightarrow \,$(1): Let $(X,\leq,\T)$ be any strongly convex C-stable semi-qospace. 
Then, by Lemma \ref{cstable}, $(X,\T^{\leq})$ is a C-space with specialization order $\leq$, 
while $\T^{\geq}$ has the dual specialization order. By strong convexity, $\T$ is the patch topology $\T^{\leq} \vee \T^{\geq}$.

By Proposition \ref{pat}, any patch space $(X,\leq, \T )$ of a web space, {\em a fortiori} of a C-space $(X,\S )$, is a web quasi-ordered space, hence ${\uparrow}$-stable.
For $x\in U \cap V$ with $U\in \S$ and $V\in \T^{\geq}$, there are $u\in U$ and $W\in \S$ with 
$x\in W \subseteq {\uparrow\!u}$. It follows that $x\in W \mathop{\cap} V \subseteq {\uparrow\!u} \mathop{\cap} V \subseteq U \mathop{\cap} V$.
Thus,  $(X,\leq, \T )$ is a sector space.

Characterization (3) of sector spaces follows from (2) and Lemma \ref{cstable}. 

\noindent The claims involving coselections $\zeta$ are now easily derived from the previous equivalences and Theorem \ref{patzeta}, 
because C-spaces are web spaces.
\EP

We are ready to establish a categorical equivalence between C-spaces and $\zeta$-sector spaces. 
As explained in \cite{Epatch}, the right choice of morphisms is a bit delicate.
Continuous maps would be the obvious morphisms between C-spaces. 
On the other hand, one would like to have as morphisms between quasi-ordered spaces the isotone continuous maps. 
But, as simple examples in \cite{Epatch} show, a continuous map between two C-spaces need not be continuous as a map between the associated $\zeta$-sector spaces.
Therefore, we take $\zeta$-proper maps (see Section \ref{convex}) as morphisms between C-spaces in order to deduce the desired isomorphism from Theorem \ref{patzeta}: 

\BT
\label{sectorcat}
For any coselection $\zeta$, the patch functor ${\rm P}_{\zeta}$ induces a concrete iso\-mor\-phism between 
the category of C-spaces and continuous (resp.\ $\zeta$-proper) maps and the category of $\zeta$-sector spaces and lower semicontinuous (resp.\ continuous) isotone maps. 
The inverse isomorphism is induced by the functor ${\rm U}$, sending any sector space $T$ to the upper C-space ${\rm U}T$.
\ET

A related class of morphisms is formed by the {\em core continuous maps}, 
i.e.\ continuous maps for which preimages of cores are cores. 
In terms of the specialization orders, the latter condition means that these maps are residual, i.e.\ upper adjoint (cf.\ \cite{Eclo} and \cite[0--1]{CLD}). 
Core continuous maps are $\alpha$-, $\sigma$- and $\upsilon$-proper, since residuated (i.e.\ lower adjoint) maps are $\alpha$-, $\sigma$- and $\upsilon$-continuous.
Topologically, a map is residuated iff preimages of point closures are point closures. From \cite{EABC} we cite:

\BCO
\label{qopen}
Via adjunction, the category of C-spaces with core continuous maps is dual to the category of C-spaces with quasi-open residuated maps.

\ECO

\section{Fan spaces}
\label{fanspaces}

We have seen that, by virtue of the weak patch functor ${\rm P}_{\upsilon}$, the C-spaces bijectively correspond
to the $\upsilon$-sector spaces. We shall now give some effective descriptions of these specific qospaces.
By a {\em fan} we mean a nonempty set of the form ${\uparrow\!u}\setminus {\uparrow\!F}$ for some finite $F$.
In case $(X,\leq,\T)$ is a semi-qospace, each fan is a $\upsilon$-sector. 
An ${\uparrow}$-stable semi-qospace in which each point has a neighborhood base of fans will be called a {\em fan space}. Thus, the fan spaces are the $\upsilon$-sector spaces.

\BX
\label{Ex51}
{\rm The Euclidean spaces $\Lambda (\mathbb{Q}^n) = (\Lambda \mathbb{Q})^n$ and $\Lambda (\mathbb{R}^n) = (\Lambda \mathbb{R})^n$ are metrizable fan spaces. 
In contrast to $\Lambda (\mathbb{R}^n)$, the rational space $\Lambda (\mathbb{Q}^n)$ is not locally compact, 
whereas the upper spaces $\Sigma (\mathbb{Q}^n)$ and $\Sigma (\mathbb{R}^n)$ are locally super\-compact.
But neither $\Sigma (\mathbb{Q}^n)$ nor $\Sigma (\mathbb{R}^n)$ is sober, since $\mathbb{Q}^n$ and $\mathbb{R}^n$ fail to be up-complete.
}
\EX

\BX
\label{Ex52}
{\rm An open web, a sector, and a fan in the Euclidean plane $\Lambda \,{\mathbb R}^2$.
\\[4mm]
\begin{picture}(300,85)
\put(10,0){
\begin{picture}(100,85)
\multiput(-6,79)(0,-2){22}{.}
\multiput(-6,79)(2,0){10}{.}
\multiput(-6,35)(2,0){18}{.}
\multiput(30,35)(0,-2){18}{.}
\multiput(30,-1)(2,0){22}{.}
\multiput(74,-1)(0,2){10}{.}
\multiput(48,19)(2,0){14}{.}
\multiput(48,19)(0,2){17}{.}
\multiput(48,53)(-2,0){17}{.}
\multiput(14,53)(0,2){14}{.}
\put(39.5,44){\circle*{3}}
\put(5,44){\line(1,0){34}}
\put(4.5,44){\circle{2}}
\put(4.5,45){\line(0,1){28}}
\put(4.5,74){\circle{2}}
\put(39.5,43){\line(0,-1){34}}
\put(39.5,8.5){\circle{2}}
\put(40,8.5){\line(1,0){30}}
\put(70.5,8.5){\circle{2}}

\put(-10,-15){\em open web, not sector}
\end{picture}
}

\put(130,0){
\begin{picture}(100,90)
\put(5,-15){\em sector, not fan}

\put(0,0){\line(0,1){80}}
\put(0,0){\line(1,0){87}}

\multiput(0,79)(2,-1){7}{.}
\multiput(14,72)(2,-.5){13}{.}
\multiput(39.5,65)(1.5,-1.5){9}{.}
\multiput(52.5,51)(1,-2){13}{.}
\multiput(65.5,25)(1,-1.5){9}{.}
\multiput(75,11)(1.5,-1.5){8}{.}
\end{picture}
}

\put(255,0){
\begin{picture}(100,90)
\put(30,-15){\em fan}

\put(-.4,0){\line(0,1){79}}


\put(0,0){\line(1,0){80}}

\multiput(-1,78)(2,0){6}{.}
\multiput(10.5,78)(0,-2){6}{.}
\multiput(12,67)(2,0){7}{.}
\multiput(26,67)(0,-2){7}{.}
\multiput(26,55)(2,0){8}{.}
\multiput(40,55)(0,-2){8}{.}
\multiput(40,41)(2,0){8}{.}
\multiput(54,41)(0,-2){8}{.}
\multiput(54,27)(2,0){8}{.}
\multiput(68,27)(0,-2){8}{.}
\multiput(68,12)(2,0){6}{.}
\multiput(78,11)(0,-2){6}{.}
\end{picture}
}
\end{picture}
}
\EX

\vspace{3.5ex}

From Proposition \ref{sectors} and Theorem \ref{sectorcat}, we infer for the case $\zeta = \upsilon$:

\BT
\label{fans}
The fan spaces are exactly the 
\vspace{-1ex}
\BIT
\IT[{\rm (1$\upsilon$)}] weak patch spaces of worldwide web spaces (C-spaces),
\IT[{\rm (2$\upsilon$)}] hyperconvex C-stable qospaces,
\IT[{\rm (3$\upsilon$)}] hyperconvex, upper regular, locally filtered and d-stable qospaces.  
\EIT
The category of C-spaces with continuous (resp.\ $\upsilon$-proper) maps 
is isomorphic to that of fan spaces and lower semicontinuous (resp.\ continuous) isotone maps. 
All fan spaces are regular and uniformizable, hence completely regular in {\sf ZFC = ZF+AC}.
\ET

Let us make the last statement in Theorem \ref{fans} more precise. A {\em quasi-uniformity} on a set $X$ is a filter $\U$ on $X \mathop{\times} X$ whose members contain the diagonal $\{ (x,x)\! : x \!\in\! X\}$ 
such that for each $U\in \U$ there is a $V \!\in \U$ with \mbox{$VV = \{ (x,z): \exists\, y\ (x \, V y\, V\, z)\} \subseteq U$.}
For a quasi-uniformity $\U$, the dual $\U^{-1}$ is obtained by exchanging first and second coordinate, and $\U^*$ denotes the uniformity generated by $\U \cup \U^{-1}$. 
The topology $\tau (\U)$ determined by $\U$ consists of all subsets $O$ such that for $x \in O$ there is a $\,U\!\in \U$ with $xU \subseteq O$.
For a comprehensive study of quasi-uniformities and their relationships to \mbox{(quasi-)} ordered spaces, refer to Fletcher and Lindgren \cite{FL}; see also K\"unzi \cite{Kue} and Nachbin \cite{Na}.
The following results are due to Br\"ummer, K\"unzi \cite{KB} and Lawson \cite{Lss}. 

\BL
\label{localunif}
If $(X,\S )$ is a locally compact space then there is a coarsest quasi-uniformity $\U$ with $\tau (\U) = \S$, and $\U$ is generated by the sets  
 
$\hspace{4ex} U\!\rightarrow\! V =\{ (u,v)  \in X\times X : u\in U \Rightarrow v\in V\} = (X\setminus U)\!\times\! X \, \cup \, X\!\times\!V$ 

\noindent with $\,U,V\in \S$ and $U\ll V$ in the continuous frame $\S$. Furthermore,
\BIT
\IT[{\rm (1)}] if $\B$ is a subbase for $\S$ such that $C = \bigcup\,\{ B\in \B : B\ll C\}$ for all $C\in \B$, 
then the sets $B\!\rightarrow\! C$ with $B,C\in \B$ and $B\ll C$ generate $\U$,
\IT[{\rm (2)}] $\tau (\U^{-1})$ is the cocompact topology $\tau_{\pi}\S$,
\IT[{\rm (3)}] $\tau(\U^*)$ is the patch topology $\S^{\pi}$.
\EIT
\EL

For C-spaces, these conclusions may be strengthened as follows.

\BPR
\label{superunif}
Let $(X,\S )$ be a C-space with specialization qoset $Q = (X,\leq )$ and interior relation $\Ro$.\,Then 
there is a coarsest quasi-uniformity $\U$ with $\tau (\U) \!= \S$,\,and
\BIT
\IT[{\rm (1)}] $\U$ is generated by the sets $y\Ro \!\rightarrow\! x\Ro = \{ (u,v) : y \mathrel{\Ro} u \Rightarrow  x \mathrel{\Ro} v \}$ with $x \mathrel{\Ro} y$,
\IT[{\rm (2)}] $\tau (\U^{-1})$ is the weak lower topology $\upsilon\widetilde{Q}$,
\IT[{\rm (3)}] the fan space $(X,\leq, \S^{\upsilon})$ is determined by $\,\U$, i.e., $\leq \ = \bigcap \U$ and $\S^{\upsilon} = \tau (\U^*)$.
Hence, the weak patch topology $\S^{\upsilon} = \S^{\pi}$ is uniformizable.
\EIT
\EPR

\BP
\mbox{We apply Lemma \ref{localunif} to the base $\B = \{ x\Ro : x \in X\}$ of $\O_{\Ro} = \S$ (Theorem \ref{Cspace}).}
By interpolation, $x \mathrel{\Ro} z$ implies $x \mathrel{\Ro} y \mathrel{\Ro} z$ for some $y$, and $x \mathrel{\Ro} y$ implies $y\Ro \ll x\Ro$ in $\O_{\Ro}$ 
(indeed, for $y\in x\Ro\subseteq \bigcup \Y$ with $\Y\subseteq \O_{\Ro}$, there is a $Y\in \Y$ with $y\in Y$, whence $y\Ro \subseteq Y\Ro = Y$). 
This yields the equation $x\Ro = \bigcup \,\{ y\Ro \in \B : y\Ro \ll x\Ro\}$, and then (1) follows from Lemma \ref{localunif}.

Now, Theorem \ref{Cspace}\,(6) and Lemma \ref{localunif}\,(2) give $\tau (\U^{-1}) = \tau_{\pi}\S = \upsilon\widetilde{Q}$.
Finally, $\bigcap \U$ is the specialization order of $\tau (\U)$, and  $\S^{\upsilon} = \S^{\pi} = \tau (\U^*)$, by Lemma \ref{localunif}\,(3).

It is also easy to check the equations $\S = \tau (\U)$ and $\upsilon\widetilde{Q} = \tau (\U^{-1})$ directly.
\EP

\vspace{1ex}

In \cite{FL} (cf.\ \cite{Kue}, \cite{Lss}, \cite{Na}), it is shown that in {\sf ZFC} a pospace \mbox{$T = (X,\leq,\T)$} is determined by a quasi-uniformity $\U$, 
i.e., $\leq \ = \bigcap \U$ and $\T = \tau (\U^*)$, iff it has a (greatest) ordered compactification, which in turn is equivalent to saying that
$T$ is {\em completely regularly ordered}, that is, a pospace such that for $x \!\in\! U \!\in\! \T$ there are continuous maps $f,g: X \rightarrow [\,0,1\,]$ with $f$ isotone and $g$ antitone,
$f(x) = g(x) = 1$, and $f(y) = 0$ or $g(y) = 0$ for $y\in X\setminus U$. We drop the antisymmetry and obtain:

\BCO
\label{complreg}
In {\sf ZFC}, every fan space is completely regularly quasi-ordered.
\ECO

From the concrete isomorphism between the category of C-spaces and that of fan spaces established in Theorem \ref{fans}, we finally infer (using Proposition \ref{product}) 
that a nonempty product of quasi-ordered spaces is a fan space iff all factors are fan spaces and almost all of them have a least element. 
Hence, in contrast to $\Lambda \,\mathbb{R}^n$ for $n < \omega$, the sequence spaces $\Lambda (\mathbb{R}^{\,\omega})$ and $(\Lambda \mathbb{R})^{\omega}$ (Example \ref{Ex44}) fail to be fan spaces.
However, $(\Lambda {\mathbb I})^{\kappa} = \Lambda ({\mathbb I}^{\,\kappa})$ is a fan space for ${\mathbb I} = [\,0,1\,]$ and any cardinal $\kappa$.

\section{Continuous domains as pospaces}
\label{domain}

We now are going to establish a representation of continuous domains by certain pospaces, generalizing the famous characterization of continuous
lattices as meet-continuous lattices whose Lawson topology is  Hausdorff (see \cite[III--2.11]{CLD}). 
In contrast to the situation of complete lattices, we require at most up-completeness.
For our purposes, we need a general kind of stability. 
Suppose $\zeta Q$ is a collection of upper sets in a qoset $Q$. A quasi-ordered space $(Q,\T) = (X,\leq,\T )$ or its topology is said to be {\em $\zeta$-stable} if
the interior operator $^{\circ}$ of the upper topology $\T^{\leq}$ satisfies 
$$
\textstyle{ Y^{\circ} = \bigcup\,\{ ({\uparrow\!y})^{\circ} : y\in Y\}} \ \mbox{ for all } Y \!\in \! \zeta Q. 
$$ 
A lower semi-qospace $(X,\leq,\T)$ is $\alpha$-stable iff the upper space $(X,\T^{\leq})$ is a C-space. 
In case $(X,\leq,\T)$ is ${\uparrow}$-stable, the operator $^\circ$ may be replaced by the interior operator of the original topology $\T$.
Thus, {\em C-stable} means {\em $\alpha$-stable plus ${\uparrow}$-stable}.

We denote by $\curlyvee Q$ and $\curlywedge Q$ the set of all finite unions (incl.\ $\emptyset \!= \bigcup \emptyset$), respectively intersections (incl.\,$X \!= \bigcap \emptyset$) of  principal filters, i.e., 
the unital join-, respectively meet-sub\-semi\-lattice they generate in the lattice $\alpha Q$ of all upper sets. 
By $\diamondsuit Q$ we denote the sublattice of $\alpha Q$\,generated by all principal filters and containing $\emptyset$\,\,and\,$X$.
\pagebreak
\BL
\label{stable}
Let $T = (Q,\T )$ be any quasi-ordered space.
\vspace{-.5ex}
\BIT
\IT[{\rm (1)}] $T$ is $\curlyvee$-stable if it is web quasi-ordered. 
\IT[{\rm (2)}] $T$ is $\curlywedge$-stable if $Q$ is a join-semilattice with least element. 
\IT[{\rm (3)}] $T$ is $\curlywedge$-stable iff the sets $\Ro y = \{ x : y\in ({\uparrow\!x})^{\circ}\}$ are ideals. 
\IT[{\rm (4)}] $T$ is $\diamondsuit$-stable iff it is $\curlyvee$-stable and $\curlywedge$-stable.
\IT[{\rm (5)}] $T$ is $\diamondsuit$-stable if $Q$ is a lattice with $0$ and continuous operations $y \mapsto x\wedge y$.
\IT[{\rm (6)}] $T$ is a $\curlyvee$-stable lower semi-qospace with locally hypercompact upper space 
iff the latter is a C-space with specialization qoset $Q$.
\EIT
\EL

\BP
(1) If $T$ is web quasi-ordered then the interior operator $^{\circ}$  preserves finite unions of upper sets \cite{Epatch}; 
hence, $({\uparrow\!F})^{\circ}\! = \bigcup \,\{ ({\uparrow\!y})^{\circ}\! : y\in F\} = \bigcup \,\{ ({\uparrow\!z})^{\circ}\! : z\in {\uparrow\!F}\}$ for finite $F$, 
because $y\leq z$ implies $({\uparrow\!z})^{\circ}\subseteq ({\uparrow\!y})^{\circ}$.

(2) In a join-semilattice $Q$ with least element, $\curlywedge Q$ is the set of all principal filters.

(3) Suppose $T$ is $\curlywedge$-stable. For finite $F \subseteq \Ro y$, the set 
\mbox{$F^{\uparrow\!} = \bigcap \, \{{\uparrow\! x} : x\in F\}$} is a member of $\curlywedge Q$, whence
$y\in \bigcap \, \{ ({\uparrow\! x})^{\circ}\! : x\in F\} = (F^{\uparrow})^{\circ} = 
\bigcup\,\{ ({\uparrow\! u})^{\circ}\! : u \in F^{\uparrow}\}$, i.e., $F^{\uparrow}\cap \Ro y \neq \emptyset$. 
In other words, $\Ro y$ is directed (and always a lower set). Conversely, if that is the case then,
for each finite subset $F$ of $Q$, one deduces from $y\in (F^{\uparrow})^{\circ}$ that $F$ is contained in $\Ro y$,
whence there is an upper bound $u$ of $F$ in $\Ro y$; thus, $y\in \bigcup\,\{ ({\uparrow\!u})^{\circ}\! : u \!\in\! F^{\uparrow}\}$.
The reverse inclusion $\bigcup\,\{ ({\uparrow\!u})^{\circ}\! : u \!\in\! F^{\uparrow}\} \!\subseteq\! (F^{\uparrow})^{\circ}$ is clear.

(4) Obviously, $\diamondsuit$-stable quasi-ordered spaces are $\curlyvee$- and $\curlywedge$-stable. Conversely,
suppose $T$ is $\curlyvee$-and $\curlywedge$-stable. Any $Y\in \diamondsuit Q$ is of the form
$Y = \bigcap \,\{ {\uparrow\!F_i} : i \! < \! n\}$ where each of the finitely many sets $F_i$ is finite. Hence,

$Y^{\circ} = \bigcap \,\{ ({\uparrow\!F_i})^{\circ} : i < n \} =$ \hfill { ($\curlyvee$-stability)} 

$\bigcap \,\{ \,\bigcup \,\{ ({\uparrow\!x_i})^{\circ} : x_i\in F_i \} : i < n \} =$ \hfill (set distributivity)

$\bigcup \,\{ \,\bigcap\,\{ ({\uparrow\!x_i})^{\circ}\! : i \!<\! n \} = (\,\bigcap\,\{ {\uparrow\!x_i}\! : i \!<\! n \})^{\circ} : x \in \prod_{i<n}\! F_i \} =$ \hfill { ($\curlywedge$-stability)} 

$\bigcup\,\{ ({\uparrow\!y})^{\circ}\! : y\in \bigcap\,\{{\uparrow\!x_i} : i<n\}, \,x \in \prod_{i<n}\! F_i \} =$ \hfill (set distributivity)

$\bigcup\,\{ ({\uparrow\!y})^{\circ}\! : y\in \bigcap\,\{ {\uparrow\!F_i} : i<n\} = Y \}.$

(5) If $Q$ is a lattice with least element and $\T$-continuous unary meet-operations $\wedge_x : y \mapsto x\wedge y$, then $T$ is web ordered (see \cite{Epatch});
hence, $T$ is $\curlyvee$-stable by (1). Furthermore, $T$ is $\curlywedge$-stable by (2), and finally $\diamondsuit$-stable by (4).

(6) Let $T =(X,\leq,\T )$ be a lower semi-qospace. Then $\leq$ is the specialization order of $\T^{\leq}$.
The upper space ${\rm U} T = (X,\T^{\leq})$ is a C-space iff  $U \!= \bigcup\,\{ ({\uparrow\!y})^{\circ}\! : y \!\in\! U\}$ for all $U \!\in\! \T^{\leq}$, 
which is equivalent to requiring that $T$ is $\curlyvee$-stable and satisfies the equation 
$U = \bigcup\,\{ ({\uparrow\!F})^{\circ} : F\subseteq U, \, F \mbox{ finite}\}$ for all $U\in \T^{\leq}$.
But the latter condition means that the upper space is locally hypercompact.
\EP

By an {\em mc-ordered space} we mean an ordered space such that every monotone net in the space 
or, equivalently, every directed subset of the space, regarded as a net, has a supremum to which it converges. 
It is shown in \cite{Epatch} that the strongly convex mc- and T$_1$-ordered spaces are exactly the patch spaces of monotone convergence spaces (d-spaces). 
We are prepared for the main result in this section, showing that Lawson spaces are determined by a convexity property, a monotone convergence property, stability properties of the interior operator, and a separation axiom.

\BT
\label{pospace}
For an ordered space $T$, the following conditions are equivalent:
\vspace{-.5ex}
\BIT
\IT[{\rm (1)}] $T$ is the Lawson space of a continuous domain.
\IT[{\rm (2)}] $T$ is the Lawson space of a meet-continuous dcpo, $\curlywedge$-stable and T$_2$.
\IT[{\rm (3)}] $T$ is a fan space whose upper space is sober.
\IT[{\rm (4)}] $T$ is a hyperconvex C-stable pospace whose upper space is a d-space.
\IT[{\rm (5)}] $T$ is hyperconvex, mc-ordered, $\diamondsuit$-stable, ${\uparrow}$-stable and T$_2$.
\EIT
\vspace{-.5ex}
In {\rm (5)} `$\,{\uparrow}$-stable and T$_2\!$' may be replaced by `T$_2$-ordered' or by `upper T$_3$-ordered'.
\ET

\BP
(1)$\,\Rightarrow\,$(3): Let $P$ be a continuous domain with $T = \Lambda P$. By Corollary \ref{Cspace}, its upper space $\Sigma P$ is a sober C-space,
and by Theorem \ref{fans}, $\Lambda P = {\rm P}_{\upsilon}\, \Sigma P$ is a fan space. 

(3)$\,\Rightarrow\,$(2): By Theorem \ref{fans}, $T =(P,\T )$ is regular, upper regular and T$_1$, hence T$_3$ (i.e.\ regular and T$_1$), upper T$_3$-ordered, and so T$_2$-ordered; moreover,
by Theorem \ref{fans} and Lemma \ref{cstable}, the upper space ${\rm U} T$ is a C-space.
By Corollary\,\,\ref{Cspace}, $P$ is a (meet-) continuous dcpo, $\T^{\leq}$ is its Scott topology $\sigma P$, 
and $\T = \T^{\leq \upsilon}$ is the Lawson topology $\lambda P$; and, again by Theorem \ref{fans}, $T$ is $C$-stable, {\it a fortiori} $\curlywedge$-stable. 

(3)$\,\Rightarrow\,$(4): Invoke Theorem \ref{fans} once more.

(2)$\,\Rightarrow\,$(5):  If $P$ is a meet-continuous dcpo then $\Sigma P$ is a web space, so 
its weak patch space, the Lawson space $\Lambda P$, is web ordered, hence ${\uparrow}$-stable (Proposition\,\ref{webpatch}); 
it is a hyperconvex semi-pospace, and mc-ordered since $\Sigma P$ is a d-space.
By Lemma \ref{stable}\,(1), $T$ is $\curlyvee$-stable, and if it is $\curlywedge$-stable, it is $\diamondsuit$-stable by Lemma \ref{stable}\,(4). 

(4)$\,\Rightarrow\,$(5): C-stable semi-pospaces are ${\uparrow}$-stable and $\alpha$-stable, hence $\diamondsuit$-stable.

(5)$\,\Rightarrow\,$(1): By Lemma \ref{stable}\,(3), the sets $\Ro y = \{ x: y\!\in\! ({\uparrow\! x})^{\circ}\}$ are directed.
As $T$ is mc-ordered, $\Ro y$ has a join $x \!=\! \bigvee \Ro y \leq y$. 
Assume $x \!<\! y$ and choose (using T$_2$) disjoint sets $V,W\in \T$ with $x\in V$ and $y\in W$. 
By hyperconvexity, there is a $U \in \T^{\leq}$ and a finite set $F$ with $x\in U\setminus {\uparrow\! F}\subseteq V$.
Then $y\in U\cap W\subseteq {\uparrow\! F}$ (as $x < y$ and $V\cap W = \emptyset$), 
and $U\cap W\in \T$ yields ${\uparrow\!(U\cap W)}\in \T^{\leq}$ if $T$ is ${\uparrow}$-stable, 
while in case $T$ is T$_2$-ordered we may assume $W = {\uparrow\!W}$, hence ${\uparrow\!(U\cap W)} = U \cap W \in \T ^{\leq}$.
Thus, $y\in ({\uparrow\! F})^{\circ} = \bigcup\,\{ ({\uparrow\! u})^{\circ}\! : u\in F\}$
by $\curlyvee$-stability, and so $u\in \Ro y$ for some $u\in F$, which leads to the contradiction $u\leq \bigvee \! \Ro y = x\in U\setminus {\uparrow\! F}$. 
Hence, $y$ is the directed join of\,\,$\Ro y$. Since $T =(P,\T )$ is mc-ordered, $\T^{\leq}$ is coarser than $\sigma P$. 
It follows that $\Ro y$ coincides with the way-below ideal $\ll\! y$; indeed, $x \ll y$ implies $x\mathrel{\Ro}y$, since $\Ro y$ is an ideal with join $y$; 
and $x\mathrel{\Ro} y$ implies $x \ll y$, since for directed $D$, $x \mathrel{\Ro} y = \!\bigvee\! D$ entails 
$y\in int_{\T^{\leq}} {\uparrow\!x} \subseteq int_{\sigma P} {\uparrow\!x} \subseteq x\!\!\ll$. 
Thus, $P$ is continuous with $\T^{\leq} = \sigma P$ (the base $\{ x\Ro : x\in P\}$ of $\sigma P$ is contained in $\T^{\leq}$), 
and $\T = \T^{\leq \,\upsilon} = \lambda P$ by hyperconvexity.
\EP

Since all topologies on a complete lattice $L$ are $\curlywedge$-stable, we obtain 
(using Proposition \ref{webpatch} and the fact that $L$ is meet-continuous iff its Scott space is a web space):

\BCO
\label{complete}
A complete lattice $L$ is continuous iff $\Lambda L$ is a T$_2$ web ordered space.
In that case, $\Lambda L$ is regular and upper T$_3$-ordered, and even compact if {\sf AC} holds.
\ECO

In categorical terminology, parts of Corollary \ref{contdom} and Theorem \ref{pospace} read as follows:

\BPR
\label{iso}
The category {\bf CD} of continuous domains and maps preserving directed joins is concretely isomorphic, 
by virtue of the Scott functor $\Sigma$, to the category {\bf SCT} of sober C-spaces, and, by virtue of the Lawson functor $\Lambda$,
to the category {\bf MFS} of mc-ordered fan spaces and isotone lower semicontinuous maps. 
The inverse isomorphism $\Lambda^-$ is induced by the forgetful functor $(P,\T) \mapsto P$. 
\EPR

\begin{picture}(300,63)(-80,0)

\put(82,54){{\bf CD}}
\put(35,10){{\bf SCT}}
\put(119,10){{\bf MFS}}
\put(76,49){\vector(-1,-1){24}}
\put(54,37){$\Sigma$}
\put(57,25){\vector(1,1){24}}
\put(71,29){$\Sigma^{-}$}
\put(106,49){\vector(1,-1){24}}
\put(121,37){$\Lambda$}
\put(125,25){\vector(-1,1){24}}
\put(99,29){$\Lambda^{\!-}$}
\put(67,10){\vector(1,0){46}}
\put(112,14){\vector(-1,0){46}}
\put(87,17){${\rm U}$}
\put(86,-1){${\rm P}_{\upsilon}$}
\end{picture}

\vspace{1ex}

\noindent We now provide examples showing that the five properties 

(1) {\em hyperconvex,} (2) {\em mc-ordered,} (3) {\em $\curlyvee$-stable,} (4) {\em $\curlywedge$-stable,} (5) {\em T$_2$-ordered,}

\noindent which together characterize Lawson spaces of continuous domains, are independent.

\BX
\label{Ex61}
{\rm
For ${\mathbb I} = [\,0,1\,] \subseteq {\mathbb R}$, the left half-open interval topology $\T = \upsilon \hspace{.2ex}{\mathbb I}\,^{\alpha}$ 
makes ${\mathbb I}$ an ${\uparrow}$-stable pospace with the five properties except hyperconvexity:
for no (open!) interval $[\,0,y\,]$ with $0 \!<\! y \!<\! 1$ and no upper set $U$, there is an open lower set 
$V \! = [\,0,\,z\,[\ = \mathbb{I} \setminus \!{\uparrow\!z}$ with $y\in U \cap V \subseteq [\,0,y\,]$, since any $x\in \ ]\,y,z\,[$ belongs to\,\,$U\cap V$.
}
\EX

\BX
\label{Ex62}
{\rm
All finite products of chains are continuous posets; but they are dcpos only if all nonempty subsets have joins. 
The pospaces $\Lambda \,{\mathbb Q}^n$ and $\Lambda \,{\mathbb R}^n$ have all five properties except of being mc-ordered, 
as neither ${\mathbb Q}^n$ nor ${\mathbb R}^n$ has a join. 
}
\EX

\BX
\label{Ex63}
{\rm
Subsets of the real square $[\,0,1\,]^{\,2}$, endowed with the Euclidean topology, are T$_2$ (and compact if closed). Suitably ordered, they provide instructive examples of T$_2$-ordered spaces. 
We sketch four such examples and list their properties. The vertical bars mean diagonals $\{ (r,r) : r\in [\,s,1\,]\}$, the horizontal bars codiagonals
$\{(r,1\!-\!r): r \in [\,0,t\,]\}$. In $T_2$ the non-filled point is eliminated.

\begin{picture}(400,75)(0,0)

\put(0,0){
\begin{picture}(100,80)

\put(0,31){\circle*{4}}

\put(29,31){\line(1,4){10}}
\put(20,31){\line(1,2){20}}
\put(10.5,31){\line(3,4){28}}
\put(29,31){\line(1,-4){10}}
\put(20,31){\line(1,-2){20}}
\put(10.5,31){\line(3,-4){30}}

\put(2,30.5){\line(1,0){25}}
\put(2,31){\line(1,0){25}}
\put(2,31.5){\line(1,0){25}}
\put(29,31){\circle*{4}}

\put(2,32){\line(1,1){35.5}}
\put(2,30){\line(1,-1){35.5}}

\put(39,69){\circle*{4}}
\put(39,-7){\line(0,1){75}}
\put(38.5,-7){\line(0,1){75}}
\put(39.5,-7){\line(0,1){75}}
\put(39,-7){\circle*{4}}

\put(49,29){$T_1$}
\end{picture}
}

\put(90,0){
\begin{picture}(100,80)

\put(0,31){\circle*{4}}

\put(29,31){\line(1,4){10}}
\put(20,31){\line(1,2){20}}
\put(10.5,31){\line(3,4){28}}
\put(29,31){\line(1,-2){9}}
\put(20,31){\line(1,-1){17}}
\put(10.5,31){\line(3,-2){27}}

\put(2,30.5){\line(1,0){25}}
\put(2,31){\line(1,0){25}}
\put(2,31.5){\line(1,0){25}}
\put(39,-7){\circle*{4}}
\put(39,12){\circle{4}}

\put(2,32){\line(1,1){35.5}}
\put(2,30){\line(2,-1){35.5}}

\put(39,69){\circle*{4}}
\put(39,-7){\line(0,1){17}}
\put(38.5,-7){\line(0,1){17}}
\put(39.5,-7){\line(0,1){17}}
\put(39,13.5){\line(0,1){55}}
\put(38.5,13.5){\line(0,1){55}}
\put(39.5,13.5){\line(0,1){55}}
\put(29,31){\circle*{4}}

\put(49,29){$T_2$}
\end{picture}
}

\put(180,0){
\begin{picture}(100,80)

\put(0,31){\circle*{4}}

\put(2,32){\line(1,1){35.5}}
\put(2,30){\line(1,-1){35.5}}

\put(39,69){\circle*{4}}
\put(39,-7){\line(0,1){75}}
\put(38.5,-7){\line(0,1){75}}
\put(39.5,-7){\line(0,1){75}}
\put(39,-7){\circle*{4}}

\put(49,29){$T_3$}
\end{picture}
}

\put(270,0){
\begin{picture}(100,80)

\put(0,31){\circle*{4}}

\put(29,31){\line(1,4){10}}
\put(20,31){\line(1,2){20}}
\put(10.5,31){\line(3,4){28}}

\put(2,30.5){\line(1,0){25}}
\put(2,31){\line(1,0){25}}
\put(2,31.5){\line(1,0){25}}
\put(29,31){\circle*{4}}

\put(2,32){\line(1,1){35.5}}

\put(39,69){\circle*{4}}
\put(39,30){\line(0,1){38}}
\put(38.5,30){\line(0,1){38}}
\put(39.5,30){\line(0,1){38}}
\put(39,31){\circle*{4}}

\put(49,29){$T_4$}
\end{picture}
}

\end{picture}

\vspace{3ex}

\noindent 
\begin{tabular}{c|c|c|c|c|c|c|c}
$T$ &$\!$hyperconvex$\!$&$\!$mc-ordered$\!$&$\!\curlyvee$-stable & $\!\curlywedge$-stable & ${\!\uparrow}$-stable & compact & complete \\
\hline
1 &$-$& + & + & + &$-$& + & + \\
2 & + &$-$& + & + &$-$&$-$&$-$\\
3 & + & + &$-$& + &$-$& + & + \\
4 & + & + & + &$-$&$-$& + &$-$
\end{tabular}

\bigskip
Notice that each point of the pospace $T_3$ has a neighborhood base of fans,
but $T_3$ fails to be a fan space because ${\uparrow}$-stability is violated (cf.\ Example \ref{Ex43}).
}
\EX

\BX
\label{Ex64}
{\rm
Let $(X,\T )$ be a nonempty T$_2$-space (like $\Lambda \,{\mathbb R}$ or $\Lambda \,{\mathbb I}$) whose finite subsets have empty interior. 
Then $(X, =, \T )$ is T$_2$-ordered, hyperconvex,  mc-ordered, $\uparrow$- and $\curlyvee$-stable, and it is almost $\curlywedge$-stable: for finite subsets $F$ with at least two elements, 
one has $F^{\uparrow} = \emptyset$, hence $(F^{\uparrow})^{\circ} = \bigcup \,\{ ({\uparrow\!y})^{\circ} :y\in F^{\uparrow}\} = \emptyset$.
But, ``by the skin of its teeth'', such a pospace fails to be $\curlywedge$-stable, because for $F = \emptyset$, the set $F^{\uparrow}$ is the whole set $X$,
whence $(F^{\uparrow})^{\circ} = X$ differs from $\bigcup \,\{ ({\uparrow\!y})^{\circ} :y\in F^{\uparrow}\! = X\} = \emptyset$.
}
\EX

\BX
\label{Ex65}
{\rm
The completed open real square $S \!= \,]\,0,1\,[^{\,2} \,\mathop{\cup}\, \{ (0,0),(1,1) \}$ is {\em not} meet-continuous.
Endowed with the interval topology, it becomes an ${\uparrow}$-stable compact semi-pospace, but clearly not T$_2$, being irreducible: 
indeed, each nonempty open set contains points close to $(1,0)$ or to $(0,1)$ (cf.\ Example 3.5 in \cite{Epatch}). 
However, $S$ is easily seen to be hyperconvex, mc-ordered, $\curlyvee$- and $\curlywedge$-stable, and so $\diamondsuit$-stable. 
}
\EX

Some of the basic results about C-spaces and fan spaces are more than 30 years old; they have been reported by the author
at the Annual Meeting of the\,\,DMV, Dortmund 1980, and at the Summer School on Ordered Sets and Universal Algebra, Donovaly 1985,
but have not been published systematically until now.

\section{Semitopological and topological semilattices}
\label{semi}

Let us now apply some of the previous results to the situation of semilattices, by which we always mean {\em meet-semilattices}. 

A {\em semilattice-ordered space} is an ordered space whose underlying poset is a semilattice;
it is said to be {\em compatible} if its specialization order is the order of the semilattice (this differs from the compatibility in \cite[VI--1]{CLD}). 
A semilattice- and T$_1$-ordered space will be referred to as a T$_1${\em -semilattice} for short. 
A {\em semitopological semilattice} is a semilattice-ordered space whose unary meet operations \mbox{$\wedge_x \!: y \mapsto x\wedge y$} are continuous,
while in a {\em topological semilattice} the binary meet operation is continuous.
Note the following equivalences for a semilattice-ordered space $S$ (see \cite{Epatch}):

{\em $S$ is semitopological $\Leftrightarrow $ $S$ is web ordered,}

\noindent and if $S$ is locally convex,

\noindent {\em $S$ has small\,(convex)\,semilattices $\,\Leftrightarrow S$ is locally filtered $\,\Leftrightarrow S$ is locally lower bounded.}

\BL
\label{top}
A semilattice-ordered space $S =(X,\leq,\T)$ is a topological semilattice whenever one of the following conditions is fulfilled:
\vspace{-.5ex}
\BIT
\IT[{\rm (1)}] $S$ is compatible and locally filtered.
\IT[{\rm (2)}] $S$ has a subbase of complements of filters (i.e.\ down-directed upper sets).
\IT[{\rm (3)}] $S$ carries the weak lower topology: $\T = \upsilon (X,\geq)$.
\IT[{\rm (4)}] $S$ is a hyperconvex, locally filtered semi-pospace.
\EIT
\EL

\BP
(1) For $x,y\in \! X$ and a filtered neighborhood $D$ of $x\mathop{\wedge} y$ with \mbox{$D \subseteq U \!\in \T \!= \! \T^{\leq}$,} the up-closure $F = {\uparrow\!D}$ is a filter still contained in $U$.
For $W = int_{\T} F$, one obtains $W\wedge W \subseteq F\wedge F \subseteq F \subseteq U$ and $x,y\in W$, as $x\wedge y \in W = {\uparrow\!W}$ by ${\uparrow}$-stability.

(2) Let $V$ be a subbasic open set such that $F = X\setminus V$ is a filter. If $x,y\in X$ satisfy $x\wedge y\in V$ then $x\in V$ or $y\in V$ (otherwise $x\wedge y \in F$).
Hence, $(x,y)$ lies in $W = (V \!\times\! X)\cup (X \!\times\! V)$, and that is an open set in $S\!\times\!S$ with $\wedge\, [W] \subseteq V$; 
indeed, $(u,v)\in W$ implies $u\in V$ or $v\in V$, and so $u\wedge v\in V$, since $V$ is a lower set. 

(3) follows from (2), because $\upsilon (X,\geq )$ has a closed subbase of principal filters.

(4) The compatible ordered space $(X,\leq,\T^{\leq})$ is locally filtered: for \mbox{$x \in U \in \T^{\leq}$,}
find a filtered $D \subseteq U$ with $x\in \!W \! =\! int_{\T}D$; then ${\uparrow\!\!D}$ is a filter with \mbox{$x \!\in\! {\uparrow\!W} \!\subseteq\! {\uparrow\!\!D} \subseteq U$,}
and ${\uparrow\!W}$ is $\T^{\leq}$-open by ${\uparrow}$-stability.
By (1), the operation $\wedge$ is $\T^{\leq}$-continuous; by (3), it is $\upsilon (X,\geq)$-continuous, 
and then, by hyperconvexity, it is $\T$-continuous.
\EP

By Lemma \ref{pat} for $\zeta = \upsilon$, the hyperconvex T$_1$-semilattices are exactly the weak patch spaces of T$_0$ spaces whose specialization poset is a semilattice.

We call a topological semilattice {\em s-topological} if it has small semilattices. Thus, a hyperconvex s-topological semilattice has even small convex semilattices.
We are ready for the characterization of hyperconvex semitopological or s-topo\-logical T$_1$-semilattices as certain web ordered spaces.

\BT
\label{topsemi}
Let $S$ be a hyperconvex T$_1$-semilattice or, equivalently, the weak patch space of a compatible semilattice-ordered space.
\begin{enumerate}[itemsep=1mm]
\item The following three conditions are equivalent:
\BIT
\IT[{\rm (11)}] $S$ is the weak patch space of a (unique) web space.
\IT[{\rm (12)}] $S$ is a web ordered space. 
\IT[{\rm (13)}] $S$ is a semitopological semilattice.
\EIT
\item The following three conditions are equivalent:
\BIT
\IT[{\rm (21)}] $S$ is the weak patch space of a (unique) wide web space.
\IT[{\rm (22)}] $S$ is a locally filtered (or locally lower bounded) ordered space. 
\IT[{\rm (23)}] $S$ is an s-topological semilattice.
\EIT
\item The following three conditions are equivalent:
\BIT
\IT[{\rm (31)}] $S$ is the weak patch space of a (unique) worldwide web space.
\IT[{\rm (32)}] $S$ is a C-stable pospace (a fan space). 
\IT[{\rm (33)}] $S$ is an s-topological semilattice whose upper space is locally compact.
\EIT
\end{enumerate}
\ET

\BP
(11)$\,\Leftrightarrow \,$(12): Apply Theorem \ref{patzeta}\,(1) to $\zeta = \upsilon$. 

(12)$\,\Rightarrow \,$(13): By Proposition \ref{webpatch}\,(1) and Theorem \ref{patzeta}\,(1), ${\rm U}S$ is a web space with $S \!=\! {\rm P}_{\upsilon}{\rm U}S$.
Therefore, the unary meet operations $\wedge_x$ are continuous in ${\rm U}S$ \cite{Eweb};
they are continuous in $\Upsilon\widetilde{S}$ by Lemma \ref{top}\,(3), and so continuous in $S$ by hyperconvexity.

(13)$\,\Rightarrow \,$(12): This was established in \cite{Epatch}.

(21)$\,\Leftrightarrow \,$(22): Apply Theorem \ref{patzeta}\,(2) to $\zeta = \upsilon$.

(22)$\,\Rightarrow \,$(23): By Lemma \ref{top}\,(4), $S = (X,\leq, \T)$ is a topological semilattice.
Given $x \!\in\! O \!\in\! \T$, use hyperconvexity and local filteredness in order to find a $U\in \T^{\leq}$, a finite set $F$, 
and a filtered set $D$ such that $x\in int_{\T} D \subseteq D \subseteq U \setminus {\uparrow\!F} \subseteq O$.
Then ${\uparrow\!D}\setminus {\uparrow\!F}$ is a (convex) subsemilattice with 
$x\in int_{\T} D \subseteq {\uparrow\!D}\setminus {\uparrow\!F} \subseteq U\setminus {\uparrow\!F} \subseteq O$,
which shows that $S$ has small (convex) semilattices.

(23)$\,\Rightarrow \,$(22): By continuity of the unary operations $\wedge_x$, $S$ is ${\uparrow}$-stable: indeed,
$O\in \T$ implies ${\uparrow\!O} = \bigcup\,\{ \wedge_x^{-1}[O] : x\in O\} \in \T$. Clearly, subsemilattices are filtered.

The equivalence of (31), (32) and (33) is easily verified with the help of (2), Theorem \ref{fans} 
and Proposition \ref{supertop}: locally compact wide web spaces are C-spaces.
\EP

\vspace{1ex}

Combining Theorems \ref{fans}, \ref{pospace} and \ref{topsemi} with Corollary \ref{contdom}, we conclude:

\BCO
\label{contsemi}
For a semilattice-ordered space $S$, the following are equivalent:
\vspace{-.5ex}
\BIT
\IT[{\rm (41)}] $S$ is the weak patch space of a (unique) sober C-space.
\IT[{\rm (42)}] $S$ is an mc-ordered fan space. 
\IT[{\rm (43)}] $S$ is a hyperconvex, mc-ordered, s-topological T$_1$-semilattice whose upper space is locally compact.
\IT[{\rm (44)}] $S$ is the Lawson space of a continuous up-complete semilattice.
\EIT
The weak patch functor ${\rm P}_{\upsilon}$ induces concrete isomorphisms between
\vspace{-.5ex}
\BIT
\IT[{\rm (1)}] the category of compatible semitopological semilattices and the category of hyper\-convex semitopological T$_1$-semilattices,
\IT[{\rm (2)}] the category of compatible s-topological semilattices and the category of hyper\-convex s-topological T$_1$-semilattices,
\IT[{\rm (3)}] the category of locally compact compatible s-topological semilattices and the category of topological semilattices that are fan spaces,
\IT[{\rm (4)}] the category of sober locally compact compatible s-topological semilattices and the category of Lawson spaces of continuous up-complete semilattices.
\EIT
\ECO

\noindent Of course, most elegant results are available for compact pospaces. From \cite[IV--1]{Com} or \cite[VI--1]{CLD}, we learn the following facts:
every compact pospace is
\vspace{-.5ex}
\BIT
\IT[(1)] mc-ordered and dually mc-ordered,
\IT[(2)] monotone normal, in particular upper regular,
\IT[(3)] strongly convex,
\IT[(4)] $\uparrow$- and $\downarrow$-costable: for any closed subset $C$, the sets ${\uparrow\!C}$ and ${\downarrow\!C}$ are closed.
\EIT
By (1), the upper space and the lower space of a compact pospace are d-spaces. But a compact pospace need not be ${\uparrow}$-stable, nor need it have small web neighborhoods (Examples \ref{Ex63},\,\ref{Ex64}).
In particular, compact pospaces need not be locally filtered.

\BPR
\label{compact}
The compact Lawson spaces of continuous domains are exactly the locally filtered compact pospaces.
\EPR

\BP
If $P$ is a continuous domain then, by Theorem \ref{pospace}, its Lawson space is a fan space, hence a locally filtered pospace, by Theorem \ref{fans}.
Conversely, if \mbox{$T = \!(X,\leq,\T)$} is a locally filtered compact pospace then $T$ is upper regular, mc-ordered and dually mc-ordered; by Lemma \ref{cstable}, it is d-stable (because ${\rm L}T$ is a d-space),
and its upper space ${\rm U}T$ is a C-space, hence a locally supercompact d-space. By Corollary \ref{contdom}, ${\rm U}T$ is the Scott space of the continuous domain $P =(X,\leq)$.
By Theorem \ref{Cspace}\,(5), the Lawson topology $\lambda P = \sigma P \vee \upsilon\widetilde{P}$ is contained in $\T^{\leq}\vee \T^{\geq}$, hence in $\T$.
Now, $\lambda P$ is T$_2$ (see Theorem \ref{pospace}) and $\T$ is compact, so both topologies coincide. 
\EP

In \cite[III--5]{CLD}, one finds a collection of various equivalent characterizations of continuous domains that are compact in their Lawson topology. 
A semilattice is said to be {\em complete} if all directed subsets have suprema and all nonempty subsets have infima (this together with the existence of a top element defines a complete lattice).
The previous results, together with the equivalence of the Ultrafilter Theorem to compactness of Lawson spaces for complete semilattices (see \cite{EPS}) immediately lead to one of the most important results of continuous semilattice theory \cite[IV--3]{CLD}:

\vspace{1ex}

\noindent {\bf Fundamental Theorem of Compact Semilattices}\\
{\em The Lawson spaces of complete continuous semilattices are exactly the compact T$_2$ s-topological semilattices.
}

\section{Domain bases and core bases}
\label{corebasis}

Following \cite{AJ} and \cite{CLD}, we mean by a {\em basis} of a poset $P$ a subset $B$ such that for each $y\in P$, the set $\{ b \in\! B: b \ll y\}$ is directed with join $y$; 
if $P$ is a domain, the pair $(P,B)$ is called a {\em based domain}. Thus, a poset is continuous iff it has a basis. 

Similarly, we call a subset $B$ of a space $(X,\S)$ a {\em core basis} if for all $U\in \S$ and all $y\in U$,
the set $\Ro_{\S}y$ meets $B\cap U$; in other words, all points have neighborhood bases formed by cores of elements of $B$.
(To avoid ambiguities, we use the word {\em basis} for subsets of the ground set and the word {\em base} for subsets of the power set.)
A space is a C-space iff it has a core basis. By a {\em core based space} we mean a pair consisting of a space and a core basis of it.
The proof of the next lemma is routine:

\BL
\label{basis}
The bases of a poset are the core bases of its Scott space. 
Hence, via the Scott functor $\Sigma$, the based domains correspond to the core based sober spaces. 
\EL

\BPR
\label{Cbasis}
The C-ordered sets are exactly the pairs $(B,\ll\!\!|_B)$, where $B$ is a basis of a domain (which is unique up to isomorphism) with way-below relation $\ll$. 
\EPR

\BP
That, for any basis $B$ of a domain, the pair $(B,\ll\!\!|_B)$ is a C-ordered set follows easily from the interpolation property (cf.\ Lemma \ref{conti}).
Conversely, sending $x$ to $\Ro x$, one obtains an isomorphism between any C-ordered set $(X,\Ro )$ and \mbox{$(\B_{\Ro}, \ll\!\!|_{\B_{\Ro}})$,}
where $\B_{\Ro}\! = \{ \Ro x \!: x \!\in\! X\}$ is a basis of the continuous domain $\I_{\,\Ro}$ of rounded ideals (see Theorem \ref{Cspace}\,(4)). 
And if $(P,B)$ is a based domain with $\Ro = \ \ll\!|_B$ then the map $y\mapsto B_y = \{ b\in B: b\ll y\}$ is an isomorphism between $P$ and $\I_{\,\Ro}$ (cf.\ \cite[2.2.6]{AJ}).
\EP

\BCO
\label{Cor}
The T$_0$ C-spaces are exactly the core bases of sober C-spaces, equipped with the induced topology. 
\ECO

Define a {\em based supercontinuous lattice} to be a pair consisting of a supercontinuous 
lattice $L$ and a {\em coprime basis}, i.e., a join-dense subset of coprime elements.
We have gathered six different descriptions of C-ordered sets and T$_0$ C-spaces:

\BT
\label{Cequiv}
The following six categories are mutually equivalent:

\vspace{1ex}

\begin{tabular}{l|l}
objects & morphisms\\
\hline
C-ordered sets & interpolating isotone maps\\
T$_0$ C-spaces & continuous maps\\
fan ordered spaces & lower semicontinuous isotone maps\\
based domains & maps preserving directed joins and bases\\
core based sober spaces & continuous maps preserving core bases\\
based supercontinuous lattices & maps preserving joins and coprime bases\\
\end{tabular}
\ET

\vspace{1ex}

\noindent On the object level, these equivalences follow from Theorems \ref{Cspace} and \ref{fans}, combined with 
Lemma \ref{Cbasis}, Proposition \ref{Cequiv}, and Corollary \ref{Cor}.
The verification of the not yet proven correspondences between the morphisms is left as an exercise.

\section{Universal constructions and powerdomains}
\label{power}

Let us add a few further categorical results on worldwide web spaces (C-spaces) that are relevant to powerdomain constructions (see \cite[6.2]{AJ} and \cite[IV-8]{CLD}).

By Proposition \ref{supertop}, for any C-space $S$, the coframe ${\rm C} S$ of all closed sets is continuous, so by Corollary \ref{contdom}, 
the Scott space $\Sigma {\rm C} S$ is a complete sober C-space (completeness refers to the specialization order). 
Let {\bf CSCT} denote the category of all complete sober C-spaces and continuous maps preserving finite joins. 
By Proposition \ref{iso}, {\bf CSCT} is isomorphic to the category {\bf CL} of continuous lattices and morphisms preserving arbitrary (i.e.\ finite and directed) joins. 
Now, the map $\eta_S : S \rightarrow \Sigma {\rm C} S, \ x \mapsto {\downarrow\!x} = \{ x\}^-$ is Scott continuous, hence continuous by Theorem \ref{Cspace}\,(5), 
and for each continuous map $f : S \rightarrow C$ into a {\bf CSCT}-object,  the map $f^{\vee} : \Sigma {\rm C} S \rightarrow C, \ A \mapsto \bigvee f[A]$ is easily seen to be the unique 
{\bf CSCT}-morphism $h$ with $f = h \circ \eta_S$ (cf.\ \cite[IV-8.5]{CLD}). Categorically speaking: 

\BPR
\label{Hoare}
The closed set lattice function ${\rm C}$ induces a reflector from the category {\bf CT} of C-spaces to the non-full subcategory {\bf CSCT}. Thus, the Scott functor $\Sigma$ from {\bf CL} to {\bf CT}
is adjoint to the functor from {\bf CT} to {\bf CL} induced by ${\rm C}$.
\EPR

Rstriction to nonempty closed sets gives the {\em lower power theory via Hoare power\-domains} \cite[IV-8.7]{CLD}. 
For the {\em upper power theory via Smyth powerdomains} \cite[IV--8.10]{CLD} both local supercompactness and sobriety are indispensible, 
so that a proper extension of the case of continuous domains is hardly available.
But we can generalize the {\em convex power theory via Plotkin powerdomains} \cite[IV--8.12]{CLD} to the setting of C-spaces. 
By a {\em C-topological semilattice} we mean a C-space equipped with a continuous semilattice operation (independent of the specialization order). 
Since the interior relation of a finite product of C-spaces is the product of the interior relations of the factors, 
the relational characterization of continuity (see Section\,\,3) for $n$-ary operations $o$ on a C-space is expressible by the equation \vspace{-.5ex}
$$\Ro\hspace{.2ex} o(y_1,...,y_n) = \Ro\hspace{.2ex} o\hspace{.2ex}[\hspace{.2ex}\Ro y_1 \!\times ...\times\! \Ro y_n].
$$
By a {\em C-semilattice} we mean a C-qoset $(X,\Ro )$ with a semilattice operation $\cdot$ satisfying
$$\Ro \hspace{.2ex}(x\cdot y) = \Ro (\Ro x \cdot \Ro y) \ (= \Ro \hspace{.2ex}\{ u\!\cdot v : u \mathrel{\Ro} x,\, v \mathrel{\Ro} y\}). \vspace{.5ex}
$$
Besides the isomorphic categories {\bf CQ} of C-qosets and {\bf CT} of C-spaces, and the isomorphic subcategories {\bf CD} of \mbox{continuous} domains (with their way-below relations) and {\bf SCT} 
of sober C-spaces, consider the following semilattice-enriched categories:

\vspace{1ex}

\noindent 
\begin{tabular}{c|l|l}
$\!$category$\!$& objects & morphisms\vspace{.5ex}\\
\hline
{\bf CQS} & C-semilattices & interpolating isotone homomorphisms$\!$\\ 
{\bf CDS} & continuous dcpo-semilattices & Scott-continuous homomorphisms\\ 
{\bf CTS} & C-topological semilattices & continuous homomorphisms\\
{\bf SCTS} & sober C-topological semilattices$\!$& continuous homomorphisms\\
\end{tabular} 

~

The semilattice operation of a dcpo-semilattice commutes with directed joins, but notice that neither binary joins nor binary meets need exist. 
The functorial concrete isomorphism ${\rm T}: {\bf CQ}\rightarrow {\bf CT}$, sending $(X,\Ro)$ to $(X,\O_{\Ro})$ (Theorem \ref{Cspace}\,(2)), induces
concrete isomorphisms ${\rm T}^{\bullet}: {\bf CQS} \rightarrow {\bf CTS}$ and $\Sigma^{\bullet}:  {\bf CDS} \rightarrow{\bf SCTS} $.

Consider a C-space $S = (X,\S)$ with interior relation $\Ro$,
or directly a C-qoset $(X,\Ro)$. On the set ${\rm F} X$ of all nonempty finite subsets of $X$, define a relation $\overline{\Ro}$ by \vspace{-1ex}
$$
E \mathrel{ \overline{\Ro}} F \ \Leftrightarrow \ E\subseteq \Ro F \mbox{ and } F \subseteq E \Ro. \vspace{-1ex}
$$
Similar arguments as in \cite[IV--8]{CLD}, with $\ll$ replaced by $\Ro$, yield the following facts:

\BL
\label{Frho}
For an arbitrary C-qoset $R =(X,\Ro )$, the relation $\overline{\Ro}$ is a C-order on ${\rm F}X$. Hence, $({\rm F} X, \overline{\Ro})$ is an abstract basis, and $({\rm F} X, \O_{\overline{\Ro}})$ is a T$_0$ C-space. 
Furthermore, ${\rm F}R = ({\rm F} X, \overline{\Ro}, \cup )$ is a C-semilattice, and $({\rm F} X, \O_{\overline{\Ro}\,}, \cup )$ is a C-topological semilattice.
\EL

The way-below functor ${\rm W}: {\bf CD} \rightarrow {\bf CQ}$, sending $(X,\leq)$ to $(X,\ll)$, is adjoint to the rounded ideal functor ${\rm I} : {\bf CQ} \rightarrow {\bf CD}$, 
sending $(X,\Ro)$ to $(\I_{\Ro},\subseteq)$, and the latter is topologically represented, via the isomorphisms ${\rm T}$ and $\Sigma$, by the sobrification reflector from {\bf CT} to {\bf SCT} (cf.\,\cite{Lri}). 
For {\bf CDS}-objects $S = (X,\leq, \cdot)$ we put ${\rm W}^{\bullet}S = (X,\ll, \cdot)$, and for {\bf CQS}-objects $S = (X,\Ro, \cdot)$ we put ${\rm I}^{\bullet} S = (\I_{\Ro}, \subseteq, \cdot_{\Ro})$, where
$I \cdot_{\Ro} J = \Ro\,\{ x\cdot y : x \in I, \,y \in J\}$.
Then we have the following categorical adjunctions:

\BT
\label{adjoint}
{\rm (1)} The sobrification gives rise to a reflector ${\rm S}^{\bullet}\!: {\bf CTS} \rightarrow {\bf SCTS}$.

{\rm (2)} ${\rm F}S = ({\rm F}X, \O_{\overline{\Ro_{\S}}},\cup)$ is the free C-topological semilattice over a C-space $S$.

{\rm (3)} The way-below functor ${\rm W}^{\bullet} : {\bf CDS} \rightarrow {\bf CQS}$ is adjoint to the rounded ideal functor ${\rm I}^{\bullet}: {\bf CQS} \rightarrow {\bf CDS}$,
with unit morphisms $\Ro_S \!: S \rightarrow {\rm W}^{\bullet}\,{\rm I}^{\bullet}S,\, x\mapsto \Ro x$.

{\rm (4)} The functor ${\rm G} : {\bf CQS} \rightarrow {\bf CQ}$ forgetting the semilattice is adjoint to the free semilattice functor ${\rm F}: {\bf CQ} \rightarrow {\bf CQS}$, 
with unit morphisms $\iota_R \!: R \rightarrow {\rm G}\hspace{.2ex}{\rm F} R, \, x \mapsto  \{x\}$.
\ET

\BP
We only prove (3) and (4). By virtue of the isomorphisms between the relational and the topological structures, (1) and (2) are then immediate (cf.\ \cite{Lri}). 

(3) ${\rm W}^{\bullet}$ is well-defined: if $S = (X,\leq,\cdot)$ is a {\bf CDS}-object then $(X, \sigma (X,\leq),\cdot)$ is a {\bf CTS}-object, and so ${\rm W}^{\bullet}S = (X,\ll,\cdot) = (X,\Ro_{\sigma(X,\leq)},\cdot)$ 
is a {\bf CQS}-object. Further,
$I\cdot_{\Ro} J = \Ro\,\{ x\cdot y : x\in \!I ,\, y\in \!J\} \in \I_{\Ro}$ for {\bf CQS}-objects $S = (X,\Ro , \cdot )$ and $I,J\in \!\I_{\Ro}$\,: given $x_i\in I$, $y_i\in J$ and
$u_i \mathrel{\Ro} (x_i\cdot y_i)$ ($i < n$), pick $x\in I$ and $y\in J$ with $x_i\leq_{\Ro} x$ and $y_i \leq_{\Ro} y$ for all $i < n$.  
Then $\Ro x_i \subseteq \Ro x$, $\Ro y_i \subseteq \Ro y$, and so $u_i \in \Ro (\Ro x_i \cdot \Ro y_i) \subseteq \Ro (\Ro x \cdot \Ro y) = \Ro (x\cdot y)$; 
now, there is a $u\in \Ro (x\cdot y) \subseteq I\cdot_{\Ro} J$ with $u_i \leq u$ for all $i$; thus, $I\cdot_{\Ro} J$ is directed.
It follows that ${\rm I}^{\bullet} S = (\I_{\Ro}, \subseteq , \cdot_{\Ro})$ is a {\bf CDS}-object: 
$(\I_{\Ro},\subseteq)$ is a continuous domain, and $\cdot_{\Ro}$ is Scott continuous, since it preserves directed joins in either coordinate; indeed, for $I\in \I_{\Ro}$ and directed $\D\subseteq \I_{\Ro}$, we compute \vspace{-.5ex}

$I \cdot_{\Ro} \bigcup \D = \Ro \,\{ x\cdot y: x\in I, \, \exists \, J\in \D \ (y\in J)\}  = \bigcup\, \{ I\cdot_{\Ro} J : J\in \D\}$.\\
$\Ro_S$ is a semilattice homomorphism: $\Ro_S (x \cdot y) = \Ro \,(\Ro x\cdot \Ro y) = \Ro_S (x) \cdot_{\Ro} \Ro_S (y)$; 
and moreover, $\Ro_S$ is a {\bf CQS}-morphism, on account of the equation $\ll \!\Ro_S (y) = \ \ll \!\Ro_S[\,\Ro y]$; 
indeed, by the equivalence $I\ll J \, \Leftrightarrow \, \exists\ x\in J \ (I\subseteq \Ro x)$ for $I,J\in \I_{\Ro}$, we have

$I\ll \Ro_S (y)\,\Leftrightarrow \exists\, x,z \,(x\mathrel {\Ro} z \mathrel{\Ro} y , I\subseteq \Ro x) \, \Leftrightarrow \exists \, z \, (z\mathrel{\Ro}y, I \!\ll \Ro_S (z)) \, \Leftrightarrow \,I\in \ \ll\!\Ro_S [\,\Ro y].$

\noindent Finally, if $f : S = (X, \Ro, \cdot) \rightarrow {\rm W}^{\bullet} T$ is a {\bf CQS}-morphism with a {\bf CDS}-object $T = (P,\odot)$, then $f^{\vee} : {\rm I}^{\bullet} S = (\I_{\Ro},\subseteq , \cdot_{\Ro}) \rightarrow T$ with
$f^{\vee} (I) = \bigvee f[I]$ is well-defined, because $I$ is directed and $f$ is isotone, whence $f[I]$ is directed, too. Furthermore,  the continuity of $f$ and the approximation property of $\ll$ yield \vspace{-.5ex}

$f^{\vee}\circ \Ro_S (y) = \bigvee f[\,\Ro y] = \bigvee \!{\ll\!f[\,\Ro y]} = \bigvee \!{\ll\!f(y)} = f(y)$, 

\noindent and by Scott continuity of $\odot$, $f^{\vee}$ is a homomorphism: \vspace{-.5ex}

$f^{\vee} (I \cdot_{\Ro} J) = \bigvee f[\,\Ro\,(I \cdot J)] = \bigvee \!\ll\! f[ I \cdot J] = \bigvee (f[I] \odot f[J]) = f^{\vee}(I)\odot f^{\vee}(J).$

\noindent Clearly, $f^{\vee}$ is the unique map $h$ that preserves directed joins and satisfies $h \circ \Ro_S = f$.

(4) For {\bf CQS}-objects $S = (X,\Ro, \cdot)$, the reduct ${\rm G} S = (X,\Ro)$ is a {\bf CQ}-object. On the other hand, for any {\bf CQ}-object $R$, Lemma \ref{Frho} assures that ${\rm F} R$ is a {\bf CQS}-object.
The map $\iota_R : R = (X,\Ro) \rightarrow {\rm G}\,{\rm F} R = ({\rm F}X, \overline{\Ro}), \, x\mapsto \{ x\}$, is a {\bf CQ}-morphism:   

$E\in \overline{\Ro} \,\iota_R (y) \, \Leftrightarrow \, E \subseteq \Ro y\, (\,\mbox{and so\,} y\in E\Ro)\, \Leftrightarrow \, \exists\, x\mathrel{\Ro} y\, (E\subseteq \Ro x) \, 
\Leftrightarrow \, E\in \overline{\Ro} \,\iota_R [\,\Ro y]$.

\noindent For any {\bf CQ}-morphism $f\! : R \!=\! (X,\Ro) \rightarrow {\rm G}S$ with $S = (X',\Ro ', \cdot)$ in {\bf CQS}, the map

$\overline{f} : {\rm F} R \rightarrow S$ with $\overline{f} (F) = f(y_1)\cdot ...\cdot f(y_n)$ for $F = \{ y_1,...,y_n\} \in {\rm F}X$

\noindent is well-defined, and $\overline{f}$ is the unique semilattice homomorphism $h$ with $h\circ \iota_R = f$. 
It remains to show that $\overline{f}$ is a {\bf CQS}-morphism, i.e., $\Ro ' \overline{f}(F) = \Ro ' \overline{f}[\,\overline{\Ro}F]$ for $F\in {\rm F}X$. \\
By continuity of $f$ and of $\cdot$, the relation $x' \mathrel{\Ro '} \overline{f}(F) = f(y_1)\cdot ...\cdot f(y_n)$ means that there exist $x_i\in X$ such that
$x_i \mathrel{\Ro} y_i$ for all $i\leq n$, and $x' \mathrel{\Ro '} (f(x_1) \cdot ... \cdot f(x_n))$. 
Putting $E = \{ x_1,...,x_n\}$, we obtain $E\mathrel{\overline{\Ro}} F$ and $x' \mathrel{\Ro '}\overline{f}[E]$.
Conversely, if there is an $E\in {\rm F}X$ with the latter two properties, find $x_i \mathrel{\Ro} y_i$ ($i\leq n$) such that $E = \{ x_1,...,x_n\}$ and $F = \{ y_1,...,y_n\}$ (repeating elements if necessary)
in order to obtain $x' \mathrel{\Ro '} \overline{f}(F)$.
\EP

\bigskip
\BC
\begin{picture}(400,162)
\put(100,110){
\begin{picture}(400,63)(0,0)
\put(82,52){{\bf CT}}
\put(5,12){{\bf SCT}}
\put(122,12){{\bf CTS}}
\put(48,-29){{\bf SCTS}}

\put(46,43){${\rm S}$}
\put(75,53){\vector(-2,-1){48}}
\put(54,36.5){$\dashv $}
\put(33,22){\vector(2,1){48}}
\put(30.5,21.3){$\iota$}
\put(64,28){${\rm E}$}

\put(89,-1){${\rm S}^{\bullet}$}
\put(120,10){\vector(-2,-1){48}}
\put(100,-6){$\dashv $}
\put(80,-20){\vector(2,1){48}}
\put(78,-20.2){$\iota$}
\put(111,-13){${\rm E}^{\bullet}$}

\put(24,-12){${\rm F}$}
\put(23,5){\vector(1,-1){24}}
\put(34,-8){$\dashv $}
\put(54,-14){\vector(-1,1){24}}
\put(44,-3){${\rm G}$}

\put(99,27){${\rm F}$}
\put(95,47){\vector(1,-1){24}}
\put(108,33){$\dashv $}
\put(127,28){\vector(-1,1){24}}
\put(119,37){${\rm G}$}

\end{picture}
}

\put(100,20){
\begin{picture}(400,63)(0,0)
\put(82,52){{\bf CQ}}
\put(5,12){{\bf CD}}
\put(122,12){{\bf CQS}}
\put(52,-29){{\bf CDS}}

\put(41,38){${\rm I}$}
\put(75,51){\vector(-2,-1){48}}
\put(46.5,31){$\dashv $}
\put(33,20){\vector(2,1){48}}
\put(56,23){${\rm W}$}

\put(94,-1){${\rm I}^{\bullet}$}
\put(122,8){\vector(-2,-1){48}}
\put(102,-8.5){$\dashv $}
\put(84,-22){\vector(2,1){48}}
\put(113,-16){${\rm W}^{\bullet}$}

\put(26,-14){${\rm P}$}
\put(23,5){\vector(1,-1){24}}
\put(34,-8){$\dashv $}
\put(54,-14){\vector(-1,1){24}}
\put(43,-3){${\rm V}$}

\put(116,37){${\rm G}$}
\put(95,45){\vector(1,-1){24}}
\put(106,33){$\dashv $}
\put(127,26){\vector(-1,1){24}}
\put(97,27){${\rm F}$}

\end{picture}
}

\put(159,28){$\Sigma^{\bullet}$}
\put(172,7){\vector(0,1){68}}
\put(174,28){$\simeq$}

\put(107,80){$\Sigma$}
\put(117,47){\vector(0,1){68}}
\put(119,80){$\simeq$}

\put(224,80){${\rm T}^{\bullet}$}
\put(237,47){\vector(0,1){68}}
\put(239,80){$\simeq$}

\put(179,126){${\rm T}$}
\put(189,86){\vector(0,1){68}}
\put(191,126){$\simeq$}

\put(270,33){\em relational}

\put(270,123){\em topological}

\end{picture}

\EC

\vspace{2.5ex}

From Theorem \ref{adjoint} one deduces by composing suitable functors (cf.\ \cite[IV--8.12]{CLD}):

\BCO
\label{Plotkin}
The forgetful functor ${\rm V}$ from {\bf CDS} to {\bf CD} is adjoint to the convex powerdomain functor 
${\rm P} = {\rm I}^{\bullet}\,{\rm F}\,{\rm W} : {\bf CD} \rightarrow {\bf CDS}, \ (X,\leq) \mapsto (\I_{\overline{\ll}} {\rm F}X, \subseteq, \cup_{\overline{\ll}})$.  
\ECO

\section{Density and weight}
\label{density}

In this last section, we are looking for characterizations of core bases in terms of the interior relation.
Recall that the lower quasi-order $\leq_{\Ro}$ of an arbitrary relation $\Ro$ on a set $X$ is given by $x \leq_{\Ro} y \,\Leftrightarrow\, \Ro x \subseteq \Ro y$.
Now, we say a subset $B$ of $X$ is 
\BIT
\IT[--] {\em $\Ro$-dense} if $x \mathrel{\Ro} y$ implies $x \mathrel{\Ro} b$ and $b \mathrel{\Ro} y$ for some $b\in B$, 
\IT[--] {\em $\Ro$-cofinal} if $x \mathrel{\Ro} y$ implies $x\!\leq_{\Ro}\!b$ and $b \mathrel{\Ro} y$ for some $b\in B$. 
\EIT
Thus, $\Ro$-cofinality means that each $B\cap \Ro y$ is cofinal in $\Ro y$ relative to $\leq_{\Ro}$. Note that 

(cd) $B$ is $\Ro$-cofinal and $\Ro$ is idempotent $\Leftrightarrow $ $B$ is $\Ro$-dense and $\Ro$ is transitive.

\vspace{1ex}

The {\em strong patch topology} $\S^{\alpha} = \S \vee \alpha (X,\geq)$ of a space $(X,\S )$ with specialization order $\leq$ coincides with the {\em Skula topology} generated by $\S \cup \S^c$ \cite{Sk}. 
Hence, the sets $C\setminus D$ with $C,D\in \S^c$ form a base for $\S^{\alpha}$.

\BPR
\label{corebase}
Let $(X,\S)$ be a C-space and $\Ro$ the corresponding interior relation.
For a subset $B$ of $X$, the following conditions are equivalent:
\BIT
\IT[{\rm (1)}] $B$ is $\Ro$-dense.
\IT[{\rm (2)}] $B$ is $\Ro$-cofinal.
\IT[{\rm (3)}] $B$ is a core basis for $(X,\S)$.
\IT[{\rm (4)}] $B$ is dense in $(X,\S^{\alpha})$ and so in every patch space of $(X,\S )$. 
\IT[{\rm (5)}] $\{ {\downarrow\!b} : b\in B\}$ is join-dense in the coframe $\S^c$ of closed sets.
\EIT
\EPR

\BP
(1)$\,\Leftrightarrow\,$(2): As $\Ro$ is idempotent by Theorem \ref{Cspace}\,(1), the equivalence (cd) applies. 

(1)$\,\Rightarrow\,$(3): $y\in U\in \S$ means $U = U\Ro$ and $x\mathrel{\Ro}y$ for some $x\in U$. 
Choose $b\in B$ with $x\mathrel{\Ro} b \mathrel{\Ro} y$.
Then $b\in \Ro y \cap U \cap B \neq \emptyset$.

(3)$\,\Rightarrow\,$(4): If $C,D$ are $\S$-closed sets with $C \!\not\subseteq\! D$, pick an $x\in C\setminus D$ 
and find a $b\in B$ with $x\in int_{\S}({\uparrow\!b}) \subseteq {\uparrow\!b} \subseteq X\setminus D$. It follows that $b\in C$,
since $b\leq x \in C = {\downarrow\!C}$, and $b\not\in D$. This proves density of $B$ in $\S^{\alpha} = \S \vee \S^c$. 

(4)$\,\Leftrightarrow\,$(5) holds for arbitrary spaces. Recall that ${\downarrow\!x}$ is the closure of the singleton $\{ x\}$ in $(X,\S)$. 
Now, \mbox{$\{ {\downarrow\!b} : b\in B\}$} is join-dense in $\S^c$
iff for $C,D\in \S^c$ with $C\not\subseteq D$ there is a $b\in B$ with ${\downarrow\!b}\subseteq C$ 
but not ${\downarrow\!b}\subseteq D$, i.e.\ $b\in C\setminus D$ (as $C$ and $D$ are lower sets).
But the latter means that $B$ meets every nonempty open set in $\S^{\alpha} = \S \vee \S^c$.

(4)$\,\Rightarrow\,$(1): Suppose $x\mathrel{\Ro} y$ and choose a $z$ with $x\mathrel{\Ro} z \mathrel{\Ro} y$. 
Then $x\Ro \in \O_{\Ro} = \S$ and therefore $z \in x\Ro \cap {\downarrow\!z} \in \S^{\alpha}$. 
Hence, there is a $b \in B \cap x\Ro \cap {\downarrow\!z}$,
and it follows that $x \mathrel{\Ro} b \mathrel{\Ro} y$, since $\Ro y$ is a lower set containing $z$ and so $b$. 
\EP

\vspace{1ex}

From now on, we assume the validity of the Axiom of Choice. Hence, each set $X$ has a cardinality, 
represented by the smallest ordinal number equipotent to\,\,$X$. 
The {\em weight} $w(\S)$ resp.\ {\em density} $d(\S)$ of a space or its topology $\S$ is the least possible cardinality of bases resp.\ dense subsets.
The weight of a C-space is at most the cardinality of any core basis $B$, since $B$ gives rise to a base $\{ b\Ro = int_{\S}({\uparrow\!b}) : b \in B\}$. 

\BL
\label{cardbase}
Every core basis of a C-space $(X,\S)$ contains a core basis of cardinality $w(\S)$. Hence, 
$w(\S)$ is the minimal cardinality of core bases for $(X,\S)$.
\EL

\BP
Let $\B$ be an arbitrary base and $B$ a core basis for $(X,\S)$. The Axiom of Choice 
gives a function picking an element $b_{\,U,V}$\,from each nonempty set of the form
$$
B_{\,U,V} = \{ b \in B : U \subseteq {\uparrow\!b} \subseteq V \} \ \ (U,V\in \B ).
$$
Then
$
B_0 = \{ b_{\,U,V} : U,V\in \B, \ B_{\,U,V} \not = \emptyset\}
$
is a subset of $B$ and still a core basis; in fact, $x\in V \in \B$ implies
\mbox{$x\in U \subseteq {\uparrow\!b} \subseteq V$} for suitable $b\in B$ and $U\in \B$, and it follows that 
$x \in U \subseteq {\uparrow\!b_{\,U,V}} \subseteq V$.
If $\S$ is infinite then so is $\B$, whence $|B_0| \leq |\B|^2 = |\B|$. Thus, we get $|B_0| \leq w(\S)$, 
and the remark before Lemma \ref{cardbase} yields equality.

If $\S$ is finite then the cores form the least base $\{ {\uparrow\!x} : x\in X\}$,
and choosing a set of representatives from these cores, one obtains a core basis of cardinality $w(\S)$.
\EP

\BX
\label{Ex101}
{\rm
Consider the ordinal space $C = \omega \!+\! 2 = \{ 0, 1, ..., \omega, \omega\!+\!1\}$ with the upper (Scott) topology.
While $\omega \cup \{ \omega\!+\!1\} = C\setminus \{ \omega \}$ is a core basis, the set $B = \omega\!+\!1 = \omega \cup \{ \omega\}$ is {\em not} a core basis:
$U \!= \{ \omega\!+\!1\}$ is $\upsilon$-open but disjoint from $B$.
However, $ \{ b\Ro \!: b \in B\} = \upsilon C \setminus \{ \emptyset\} $ is a base; note $x\Ro ={\uparrow\!x}$ for $x\neq  \omega$, but $\omega\Ro = \{ \omega\!+\!1\}$.
}
\EX
Generalizing the weight of topologies, one defines the {\em weight} $w(P)$ of a poset $P$ as the least possible cardinality of join-dense subsets of $P$.
For any relation $\Ro$ on a set $X$, we define the {\em $\Ro$-cofinality}, denoted by $c(X,\Ro)$, 
to be the minimal cardinality of $\Ro$-cofinal subsets of $X$. If $\Ro$ is idempotent then $c(X,\Ro)$
is also the {\em $\Ro$-density}, the least cardinality of $\Ro$-dense subsets.
From Proposition \ref{corebase} and Lemma \ref{cardbase} we infer:

\BT
\label{w}
For a C-space $(X,\S)$ with interior relation $\Ro$,
\vspace{-.5ex}
$$
c(X,\Ro) = w(\S) = w(\S^c) = w(\S^{\upsilon}) = d(\S^{\alpha}). 
$$
\ET

\noindent Indeed, if $B$ is a core basis for $(X,\S)$ then $\B = \{ b\Ro \setminus {\uparrow\!F} : b\in B,\, F \!\subseteq\! B, \, F \,\mbox{finite}\}$ is a base for $\S^{\upsilon}$ 
(because for $x\not\leq y$ there is a $b\in B$ with $b\leq x$ but $b\not\leq y$),
and if $B$ is infinite, $\B$ and $B$ have the same cardinality. If $B$ is finite, the base $\{ {\uparrow\!b} : b \in B\}$ of $\S$ 
is equipotent to the base $\{ {\uparrow\!b}\cap {\downarrow\!b} : b \in B\}$ of $\S^{\upsilon}.$ 

Since in {\sf ZCF} any supercontinuous lattice is isomorphic to the topology of a C-space \cite{EABC}, 
Theorem \ref{w} entails a fact that was shown choice-freely in \cite{ECD}:

\BCO
The weight of a supercontinuous (i.e.\ completely distributive) lattice is equal to the weight of the dual lattice (which is supercontinuous, too). 
\ECO

\BX
\label{Ex102}
{\rm
On the real line ${\mathbb R}$ with the upper (Scott) topology $\S$, the interior relation is the usual $<\,$. 
The rationals form a core basis, being $<$-dense in ${\mathbb R}$. Hence,
$$\omega = c({\mathbb R},<) = w(\S) = w(\S^c) = w(\S^{\upsilon}) = d(\S^{\alpha}) < w(\S^{\alpha})
$$
because the half-open interval topology $\S^{\alpha}$ has no countable base. Furthermore,
$$\omega = d(\alpha {\mathbb R}\,^\upsilon) = d(\S^{\alpha})< d(\alpha {\mathbb R}\,^\alpha ) = w(\alpha {\mathbb R}\,^\alpha ) = w(\P\, {\mathbb R}) =  |{\mathbb R}|.$$
}
\EX

\end{document}